\documentclass[final,1p,times]{article}
\usepackage[utf8]{inputenc}
\usepackage[T1]{fontenc}

\usepackage[a4paper,
            left=1in,
            right=1in]{geometry}
\usepackage{amsmath}
\usepackage{graphicx}
\usepackage{tabularx}
\usepackage[colorlinks=true, allcolors=blue]{hyperref}
\usepackage{pgfplots}
\usepackage{amsmath}
\usepackage{amssymb}
\usepackage{amsthm}
\usepackage{mathtools}
\usepackage{bigints}
\usepackage{wasysym}
\usepackage{graphicx}
\usepackage{caption}
\usepackage{subcaption}
\usepackage[noabbrev,nameinlink]{cleveref}
\crefname{algocf}{algorithm}{algorithms}
\Crefname{algocf}{Algorithm}{Algorithms}
\Crefname{remark}{Remark}{Remarks} 
\Crefname{table}{Table}{Tables}
\Crefname{enumi}{Step}{Steps}
\crefname{enumi}{step}{steps}
\usepackage{booktabs}
\usepackage{todonotes}
\usepackage{bm}
\usepackage{bbm}
\usepackage{color}
\usepackage{setspace}
\usepackage{comment}
\usepackage{multirow}
\usepackage{longtable}
\usepackage{enumerate}
\usepackage[shortlabels]{enumitem}
\usepackage{pgfplots}
\pgfplotsset{compat=1.18}
\usepackage{pgfplotstable}
\usepackage[ruled,vlined,linesnumbered,algo2e]{algorithm2e}
\usepackage{tikz}
\usepackage{nicefrac}
\usepackage{siunitx}

\definecolor{pointCloud}{rgb}{0, 0.4470, 0.7410}
\definecolor{ourSolution}{RGB}{126, 47, 142}
\definecolor{theirSolution}{RGB}{217, 83, 25}

\usepgfplotslibrary{groupplots} 
\usetikzlibrary{pgfplots.groupplots} 
\usetikzlibrary{spy}

\usepgfplotslibrary{colorbrewer}

\newif\ifblackandwhitecycle
\gdef\patternnumber{0}

\pgfkeys{/tikz/.cd,
	zoombox paths/.style={
		draw=orange,
		very thick
	},
	black and white/.is choice,
	black and white/.default=static,
	black and white/static/.style={ 
		draw=white,   
		zoombox paths/.append style={
			draw=white,
			postaction={
				draw=black,
				loosely dashed
			}
		}
	},
	black and white/static/.code={
		\gdef\patternnumber{1}
	},
	black and white/cycle/.code={
		\blackandwhitecycletrue
		\gdef\patternnumber{1}
	},
	black and white pattern/.is choice,
	black and white pattern/0/.style={},
	black and white pattern/1/.style={    
		draw=white,
		postaction={
			draw=black,
			dash pattern=on 2pt off 2pt
		}
	},
	black and white pattern/2/.style={    
		draw=white,
		postaction={
			draw=black,
			dash pattern=on 4pt off 4pt
		}
	},
	black and white pattern/3/.style={    
		draw=white,
		postaction={
			draw=black,
			dash pattern=on 4pt off 4pt on 1pt off 4pt
		}
	},
	black and white pattern/4/.style={    
		draw=white,
		postaction={
			draw=black,
			dash pattern=on 4pt off 2pt on 2 pt off 2pt on 2 pt off 2pt
		}
	},
	zoomboxarray inner gap/.initial=5pt,
	zoomboxarray columns/.initial=2,
	zoomboxarray rows/.initial=2,
	subfigurename/.initial={},
	figurename/.initial={zoombox},
	zoomboxarray/.style={
		execute at begin picture={
			\begin{scope}[
				spy using outlines={%
					zoombox paths,
					width=\imagewidth / \pgfkeysvalueof{/tikz/zoomboxarray columns} - (\pgfkeysvalueof{/tikz/zoomboxarray columns} - 1) / \pgfkeysvalueof{/tikz/zoomboxarray columns} * \pgfkeysvalueof{/tikz/zoomboxarray inner gap} -\pgflinewidth,
					height=\imageheight / \pgfkeysvalueof{/tikz/zoomboxarray rows} - (\pgfkeysvalueof{/tikz/zoomboxarray rows} - 1) / \pgfkeysvalueof{/tikz/zoomboxarray rows} * \pgfkeysvalueof{/tikz/zoomboxarray inner gap}-\pgflinewidth,
					magnification=3,
					every spy on node/.style={
						zoombox paths
					},
					every spy in node/.style={
						zoombox paths
					}
				}
				]
			},
			execute at end picture={
			\end{scope}
			\node at (image.north) [anchor=north,inner sep=0pt] {\subcaptionbox{}{\phantomimage}};
			\node at (zoomboxes container.north) [anchor=north,inner sep=0pt] {\subcaptionbox{}{\phantomimage}};
			\gdef\patternnumber{0}
		},
		spymargin/.initial=0.5em,
		zoomboxes xshift/.initial=1,
		zoomboxes right/.code=\pgfkeys{/tikz/zoomboxes xshift=1},
		zoomboxes left/.code=\pgfkeys{/tikz/zoomboxes xshift=-1},
		zoomboxes yshift/.initial=0,
		zoomboxes above/.code={
			\pgfkeys{/tikz/zoomboxes yshift=1},
			\pgfkeys{/tikz/zoomboxes xshift=0}
		},
		zoomboxes below/.code={
			\pgfkeys{/tikz/zoomboxes yshift=-1},
			\pgfkeys{/tikz/zoomboxes xshift=0}
		},
		caption margin/.initial=4ex,
	},
	adjust caption spacing/.code={},
	image container/.style={
		inner sep=0pt,
		at=(image.north),
		anchor=north,
		adjust caption spacing
	},
	zoomboxes container/.style={
		inner sep=0pt,
		at=(image.north),
		anchor=north,
		name=zoomboxes container,
		xshift=\pgfkeysvalueof{/tikz/zoomboxes xshift}*(\imagewidth+\pgfkeysvalueof{/tikz/spymargin}),
		yshift=\pgfkeysvalueof{/tikz/zoomboxes yshift}*(\imageheight+\pgfkeysvalueof{/tikz/spymargin}+\pgfkeysvalueof{/tikz/caption margin}),
		adjust caption spacing
	},
	calculate dimensions/.code={
		\pgfpointdiff{\pgfpointanchor{image}{south west} }{\pgfpointanchor{image}{north east} }
		\pgfgetlastxy{\imagewidth}{\imageheight}
		\global\let\imagewidth=\imagewidth
		\global\let\imageheight=\imageheight
		\gdef\columncount{1}
		\gdef\rowcount{1}
		
	},
	image node/.style={
		inner sep=0pt,
		name=image,
		anchor=south west,
		append after command={
			[calculate dimensions]
			node [image container,subfigurename=\pgfkeysvalueof{/tikz/figurename}-image] {\phantomimage}
			node [zoomboxes container,subfigurename=\pgfkeysvalueof{/tikz/figurename}-zoom] {\phantomimage}
		}
	},
	color code/.style={
		zoombox paths/.append style={draw=#1}
	},
	connect zoomboxes/.style={
		spy connection path={\draw[draw=none,zoombox paths] (tikzspyonnode) -- (tikzspyinnode);}
	},
	help grid code/.code={
		\begin{scope}[
			x={(image.south east)},
			y={(image.north west)},
			font=\footnotesize,
			help lines,
			overlay
			]
			\foreach \x in {0,1,...,9} { 
				\draw(\x/10,0) -- (\x/10,1);
				\node [anchor=north] at (\x/10,0) {0.\x};
			}
			\foreach \y in {0,1,...,9} {
				\draw(0,\y/10) -- (1,\y/10);
				\node [anchor=east] at (0,\y/10) {0.\y};
			}
		\end{scope}    
	},
	help grid/.style={
		append after command={
			[help grid code]
		}
	},
}

\newcommand\phantomimage{%
	\phantom{%
		\rule{\imagewidth}{\imageheight}%
	}%
}
\newcommand\zoombox[2][]{
	\begin{scope}[zoombox paths]
		\pgfmathsetmacro\xpos{
			(\columncount-1)*(\imagewidth / \pgfkeysvalueof{/tikz/zoomboxarray columns} + \pgfkeysvalueof{/tikz/zoomboxarray inner gap} / \pgfkeysvalueof{/tikz/zoomboxarray columns} ) + \pgflinewidth
		}
		\pgfmathsetmacro\ypos{
			(\rowcount-1)*( \imageheight / \pgfkeysvalueof{/tikz/zoomboxarray rows} + \pgfkeysvalueof{/tikz/zoomboxarray inner gap} / \pgfkeysvalueof{/tikz/zoomboxarray rows} ) + 0.5*\pgflinewidth
		}
		\edef\dospy{\noexpand\spy [
			#1,
			zoombox paths/.append style={
				black and white pattern=\patternnumber
			},
			every spy on node/.append style={#1},
			x=\imagewidth,
			y=\imageheight
		] on (#2) in node [anchor=north west] at ($(zoomboxes container.north west)+(\xpos pt,-\ypos pt)$);}
		\dospy
		\pgfmathtruncatemacro\pgfmathresult{ifthenelse(\columncount==\pgfkeysvalueof{/tikz/zoomboxarray columns},\rowcount+1,\rowcount)}
		\global\let\rowcount=\pgfmathresult
		\pgfmathtruncatemacro\pgfmathresult{ifthenelse(\columncount==\pgfkeysvalueof{/tikz/zoomboxarray columns},1,\columncount+1)}
		\global\let\columncount=\pgfmathresult
		\ifblackandwhitecycle
			\pgfmathtruncatemacro{\newpatternnumber}{\patternnumber+1}
			\global\edef\patternnumber{\newpatternnumber}
		\fi
	\end{scope}
}

\newcommand{\footremember}[2]{%
    \footnote{#2}
    \newcounter{#1}
    \setcounter{#1}{\value{footnote}}%
}
\newcommand{\footrecall}[1]{%
    \footnotemark[\value{#1}]%
}

\def\calP{\mathcal{P}}
\def\Pnws{\mathcal{P}_{n}^{\star}}

\theoremstyle{definition}

\newtheorem{remark}{Remark}

\DeclareMathOperator{\spn}{span}
\DeclareMathOperator{\argmin}{argmin}

\AtBeginDocument{}

\theoremstyle{plain}
\newtheorem{theorem}{Theorem}
\newtheorem{corollary}{Corollary}

\date{}

\title{A general formulation of reweighted least squares fitting}
	\author{Carlotta Giannelli\footremember{a}{Dipartimento di Matematica e Informatica, Universit{\`a} degli Studi di Firenze, \texttt{carlotta.giannelli@unifi.it}, \texttt{sofia.imperatore@unifi.it}}
	\and
	Sofia Imperatore\footrecall{a}
	\and
	Lisa Maria Kreusser\footremember{b}{Department of Mathematical Sciences, University of Bath, \texttt{lmk54@bath.ac.uk}}
	\and
	Estefan{\'i}a Loayza-Romero\footremember{c}{Department of Mathematics, Imperial College London, \texttt{k.loayza-romero@imperial.ac.uk}}
	\and
	Fatemeh Mohammadi\footremember{d}{Departments of Mathematics and Computer Science, KU Leuven, \texttt{fatemeh.mohammadi@kuleuven.be}}
	\and
	Nelly Villamizar\footremember{e}{Departments of Mathematics, Swansea University, \texttt{n.y.villamizar@swansea.ac.uk}}%
	}

\begin{document}

\maketitle

\begin{abstract}
We present a generalized formulation for reweighted least squares approximations. The goal of this article is twofold: firstly, to prove that the solution of such problem can be expressed as a convex combination of certain interpolants when the solution is sought in any finite-dimensional vector space; secondly, to provide a general strategy to iteratively update the weights according to the approximation error and apply it to the spline fitting problem. In the experiments, we provide numerical examples for the case of polynomials and splines spaces. Subsequently, we evaluate the performance of our fitting scheme for spline curve and surface approximation, including adaptive spline constructions.
\end{abstract}
	

\section{Introduction}
\label{sec:intro}
We consider the problem of constructing a continuous model $\bm v\colon \Omega \to \mathbb{R}^D$ 
on a general domain  $\Omega\subseteq\mathbb{R}^N$ in any dimension $N\in \mathbb{N}_{>0}$ from pointwise data $\left\{\left(\bm x_1, \bm f_1\right), \ldots, \left(\bm x_m, \bm f_m\right)\right\}$,
where $\bm f_i \in\mathbb{R}^D$ for $D\in \mathbb{N}_{>0}$ are
observations at points $\bm x_i\in \Omega$, for $i=1,\ldots,m$.
Given the ubiquity of this problem across various fields and applications, a number of methods have been developed to address approximation, data fitting, estimation, and prediction. Examples include classical ones, like interpolation and least squares, while more recently, approaches such as compressive sensing or neural networks, have also emerged as viable options for tackling these challenges. 
In this work, we primarily concentrate on two very well established classical methods: interpolation and {wei\-ghted} least squares, see, e.g., \cite{davis1975}. Even though interpolation is well established, there is still ongoing research in this field, including the development of efficient algorithms with irregularly-spaced data, see e.g., \cite{CAVORETTO2019,Cavoretto2020}, among others. Weighted least squares methods \cite{Strutz2015} can be considered a more general instance of ordinary least squares methods, studied on weakly admissible meshes in \cite{Bos2010}. The key difference between weighted and ordinary least squares is that in the former, fixed weight values are associated with the observations to incorporate different weightings in the least squares scheme.

Most of the time, interpolation and least squares methods are regarded as complementary techniques; however, despite this perception they share a strong connection in their polynomial formulation, as noted in \cite{floater2020}.  Additionally, as detailed in \cite{DEMARCHI2015,DELLACCIO2022} they can also be used in combination to defeat the Runge phenomenon \cite{runge1901}.
In the first part of this paper, we propose a general formulation of the weighted least {squares} approximant as a convex combination of suitable interpolants, for any finite-dimensional function space consisting of real-valued functions defined on a domain $\Omega\subseteq \mathbb{R}^N$ and address its consequences. As a special case, our formulation includes the vector space of polynomials up to a certain degree, for which a relation between weighted polynomial least squares and interpolation has been discussed in \cite{floater2020}.
In the second part of the paper, we focus on the weights of the problem. Our main aim is to update these weights in the spirit of the Iterative Reweighted Least Squares (IRLS) method~\cite{osborne1985,Bjorck1996} with convergence guarantees and efficient algorithms \cite{wolke1988,beck2015}. IRLS also offers robust regression \cite{Holland1977} and is used to smooth the reconstruction. However, the main difference between our approach and the IRLS method is that our updates allow us to preserve sharp features of the final model and are also suitable for adaptive approximations. Our numerical experiments show the performance of the proposed method within the spline framework, from curve to surface fitting with an appropriate underlying mesh, including adaptive spline constructions.

The paper is organized as follows. In \Cref{sec:int-wls}, we provide a brief introduction to the concepts of interpolation and weighted least squares approximation for a given set of observations.
In \Cref{sec:theorem} we present our theoretical result and address its consequences.
In \Cref{sec:splines}, we describe the spline models together with their hierarchical extensions used later on in this paper.
Moreover, we present a numerical verification of our theoretical result for the space of polynomial splines and comment on its consequences by exploiting hierarchical box-splines.
In \Cref{section:rWLS}, we describe the reweighted least squares spline fitting scheme and thereafter its extension to deal with adaptive spline constructions.
Finally, in \Cref{sec:Numerical_Experiments}, we report three numerical experiments to evaluate the performance of the proposed reweighted least squares fitting schemes.

\section{Interpolation and weighted least squares}
\label{sec:int-wls}
Throughout this paper, we assume that $\mathbb R^N$, $\mathbb R^D$, $\mathbb R^n$ and $\mathbb R^m$ are column vectors. Let $\Omega\subseteq\mathbb{R}^N$ and consider a set of given observations $\left\{\left(\bm x_1, \bm f_1\right), \ldots, \left(\bm x_m, \bm f_m\right)\right\}$, where $\bm x_i = (x_i^1, \ldots, x_i^N)\in \Omega$ and $\bm f_i = (f_i^1, \ldots, f_i^D)\in\mathbb{R}^D$.
Moreover, let $V$ be a vector space of functions defined on $\Omega$ and taking values in $\mathbb R^D$. We denote its finite dimension by $n = \dim V$ and assume that $V$ is generated by a basis $\Gamma = \{\beta_1, \ldots, \beta_n\}$, with basis functions $\beta_j: \Omega \to \mathbb{R}^D$, i.e.\ $V = \spn \left\{\Gamma\right\}$.

The \emph{interpolation problem} aims to find an element $\bm v\in V$ such that $\bm v\left(\bm x_i\right) = \bm f_i$ holds for each $i=1, \ldots, m$.
If a solution $\bm v$ exists,  it can be expressed as $\bm v = \sum_{j=1}^n\bm c_j \beta_j$ where the coefficients $\bm c_j = (c_j^1,\ldots,c_j^D) \in \mathbb{R}^D$ can be determined by solving a linear system for each component $k=1, \ldots, D$. Specifically, each linear system takes the  form $B\bm c^k = \bm f^k$, where
\begin{equation}\label{eq:colloc}
	B=\begin{pmatrix}
		\beta_1(\bm x_1)   & \dots  & \beta_n(\bm x_1) \\
		\vdots        &   \ddots     & \vdots      \\
		\beta_1(\bm x_m) & \dots & \beta_n(\bm x_m)
	\end{pmatrix} \in \mathbb{R}^{m\times n}
\end{equation}
is the $m\times n$ collocation matrix,
$\bm c^k=( c^k_1,\ldots, c^k_n)\in \mathbb{R}^{n}$ and $\bm f^k=( f^k_1,\ldots, f^k_m)\in \mathbb{R}^{m}$. The solution $\bm v\in V$ is unique, if and only if $B$ is invertible, hence $n=m$ is a necessary condition.
However, even if the interpolant $\bm v\in V$ exists, it is well known that it can be affected by  poor approximation quality. Several factors can contribute to this, such as the nature of the data (e.g., their quantity or distribution) and the intrinsic properties of the function space $V$ (e.g., the presence of Runge's phenomenon).

As an alternative to interpolation, a common approach is to define the approximant $\bm v\in V$ as the solution to a \emph{weighted least squares problem} when $n\leq m$.
In this context,  we assign a strictly positive \emph{weight} value $\omega_i\in\mathbb{R}_{>0}$ to each $\bm x_i\in \Omega$ for $i=1,\ldots,m$. Let $W\in \mathbb{R}^{m\times m}$ be the associated diagonal matrix with the $i$-th diagonal entry given by $\omega_i$. Then, the weighted least squares problem involves finding a function $\bm v\in  V$, which solves the  {minimization problem}
\begin{equation}\label{eq:Vwls}
	\underset{\bm u\in V}{\min}\sum_{i=1}^m\omega_i\left\|\bm u(\bm x_i) -\bm f_i\right\|_2^2.
\end{equation}
By expanding the solution $\bm v \in V$ of \eqref{eq:Vwls} as $\bm v=\sum_{j=1}^n\bm c_j\beta_j$,  the matrix $\bm c=\left(\bm c_1,\ldots,\bm c_n\right)^\top \in\mathbb{R}^{n\times D}$ is the solution to  problem \eqref{eq:Vwls} in its equivalent \emph{normal form},
\begin{equation}\label{eq:ls-normalform}
	\min_{\bm c\in\mathbb{R}^{n\times D}} \| W^{\frac{1}{2}}B\bm c - W^{\frac{1}{2}}\bm f\|_2^2,
\end{equation}
{where, $\bm f =\left(\bm f_1,\ldots,\bm f_m\right)^\top \in\mathbb{R}^{m\times D}$, }
\[
	W = \textup{diag}\{\omega_1,\ldots, \omega_m\} =  \begin{pmatrix}
		\omega_1    &  0                &  \dots  & 0  \\
		0                  & \omega_2  &  \dots  & 0 \\
		\vdots          &     \vdots    &   \ddots          & \vdots\\
		0                  &  0                & \dots   & \omega_m
	\end{pmatrix}\in\mathbb{R}^{m\times m}
\] and $W^{\frac{1}{2}}$ denotes the element-wise evaluation of the square root of $W$.
The minimization in \eqref{eq:ls-normalform} entails solving a system of linear equations for each column $\bm c^k = (c_1^k, \ldots, c_n^k)$, $k=1, \ldots, D$, of $\bm c$. More precisely, for each $k=1, \ldots, D$, $\bm c^k$  is the solution of the normal equation
\begin{equation}\label{eq:ls-normaleq}
	B^\top WB\bm c^k =  B^\top  W \bm f^k.
\end{equation}

\section{{Weighted least squares via interpolation}}
\label{sec:theorem}
In this section, to alleviate the notation, we will assume that $D=1$. Therefore, we avoid the superscript $k$ introduced in \Cref{sec:int-wls} and reduce the bold writing according to the dimension. Nevertheless, all the results apply in the more general case of $D\geq1$.
Consider $\Omega\subseteq\mathbb{R}^N$,
$\left\{\left(\bm x_1,  f_1\right), \ldots, \left(\bm x_m, f_m\right)\right\}$ a set of given observations with $\bm x_i\in\Omega$ and $f_i\in\mathbb{R}$ and finally $V = \spn \left\{\Gamma\right\}$ a vector space of dimension $n$, with basis $\Gamma = \{\beta_1, \ldots, \beta_n\}$ and functions $\beta_j:\Omega\to\mathbb{R}$, for $j=1, \ldots, n$.
For any set $K\subset \mathbb N$ we denote by $\#(K)$ the cardinality of $K$ and we define $\calP_{n}=\{K\subseteq \{1,\dots, m\}\colon \#(K)=n\}$, the set of all subsets of $\{1,\dots, m\}$ with cardinality $n$.
For each $K\in\calP_{n}$, there exist  $k_i\in \{1,\ldots,m\}$ for $i=1,\ldots,n$, such that $K=\{k_1,\ldots,k_n\}$. We denote by $v_K\in V$ the interpolant of $f_{k_i}$ at  points $\bm x_{k_i}$ for $i=1,\dots, n$.
For such $K$, we define
\begin{equation}\label{eq:colloc-loc}
	B_K=\begin{pmatrix}
	\beta_1(\bm x_{k_1})   & \dots  & \beta_n(\bm x_{k_1}) \\
	\vdots        &   \ddots     & \vdots      \\
	\beta_1(\bm x_{k_m}) & \dots & \beta_n(\bm x_{k_m})
\end{pmatrix} \in \mathbb{R}^{n\times n}
\end{equation}
which is called the \emph{collocation matrix} at the points $\bm x_{k_i}$ with respect to the whole basis $\Gamma$ of $V$.
For $K\in\calP_{n}$ and positive weights $\omega_1,\dots, \omega_m\in \mathbb{R} $, we write 
\begin{equation}\label{eq:prodwk}
\omega_K=\prod_{i\in K} \omega_i.
\end{equation}
Following the notation in \cite{floater2020}, we define
$\lambda_K=\omega_K\left|B_K\right|^2$, where $\bigl|B_K\bigr|$ denotes the determinant of   $B_K$. Finally, we  define $\calP_{n}^{\star} \subset\calP_{n}$ as the set of all subsets of $\{1,\dots,m\}$ of cardinality $n$ with $|B_K|\neq 0$. In other words, $\calP_{n}^{\star}=\{K\subseteq \{1,\dots, m\}\colon \#(K)=n,\ \text{and}\ |B_K|\neq 0\}$.
\begin{theorem}\label{th:genweightedsum}
	The weighted least squares approximant $v\in V$ of the set of points $\bigl\{(\bm x_i,f_i)\bigr\}_{i=1}^m$ is the weighted sum of the interpolants $v_K\in V$ for $K\in\Pnws$ which interpolate the points $(\bm x_i,f_i)$ for $i\in K$, i.e.,
	\begin{equation}
		\label{eq:001}
		v(\bm x) = \frac{\sum_{K\in \Pnws}\lambda_Kv_K(\bm x)}{\sum_{K\in \calP_{n}^{\star}}\lambda_K},\quad \forall \bm x\in\mathbb{R}^N.
	\end{equation}
\end{theorem}
\begin{proof}
	By hypothesis, $v\in V$  minimizes \eqref{eq:Vwls}. If $\Gamma = \left\{\beta_1, \ldots, \beta_n\right\}$ is a basis of $V$, then we can write
	$v(\bm x) = \sum_{j=1}^{n}c_j\beta_j(\bm x)$,
	for some coefficients $c_j\in\mathbb{R}$ for $j=1,\ldots,n$. We set $\bm c =(c_1,\dots, c_{n})^\top$ and $\bm f=(f_1,\dots, f_m)^\top$, and write the normal equation \eqref{eq:ls-normaleq}
	as
	$A\bm{c} = \bm{b}$ with
	\[
	\begin{split}
		A & =B^\top WB=
		\begin{pmatrix}
			\omega_1\beta_1(\bm x_1)   & \dots  & \omega_m\beta_1(\bm x_m) \\
			\vdots        &   \ddots     & \vdots      \\
			\omega_1\beta_n(\bm x_1) & \dots & \omega_m\beta_n(\bm x_m)
		\end{pmatrix}
		\begin{pmatrix}
			\beta_1(\bm x_1)    & \dots  & \beta_n(\bm x_1)  \\
			\vdots & \ddots &         \vdots \\
			\beta_1(\bm x_m)    & \dots  & \beta_n(\bm x_m)
		\end{pmatrix},\\
		\bm{b} & = B^\top W\bm{f}=\begin{pmatrix}
			\sum_{i=1}^m\omega_i\beta_1(\bm x_i)f_i, & \ldots, & \sum_{i=1}^m\omega_i\beta_n(\bm x_i)f_i
		\end{pmatrix}^\top.
	\end{split}
	\]
	By the Cauchy-Binet theorem \cite[Section 4.6]{broida1989}, $\left|A\right|
	= \sum_{K\in \calP_{n}}\omega_K\left|B_K\right|^2$
	and by Cramer's rule, we obtain $c_j  = \left|A_j\right| \Big/ \left|A\right|$, for each $j=1,\dots,n$,
	where $A_j$ is obtained from $A$ replacing its $j$-th column with $\bm{b}$ for $j=1,\dots,n$, namely, $A_j = B^\top WB_j$ with
	\[
	B_j=\begin{pmatrix}
		\beta_1(\bm x_1)  &  \dots  & \beta_{j-1}(\bm x_1)  &f_1 & \beta_{j+1}(\bm x_1)& \dots   & \beta_n(\bm x_1)  \\
		\vdots &        & \vdots & \vdots&\vdots&&\vdots \\
		\beta_1(\bm x_m)    & \dots  & \beta_{j-1}(\bm x_m)  &f_m & \beta_{j+1}(\bm x_m)& \dots   & \beta_n(\bm x_m)
	\end{pmatrix}\in \mathbb{R}^{m\times n}.
	\]
	Again from the Cauchy-Binet theorem, we obtain
	$\left|A_j\right| = \sum_{K\in \calP_{n}}\omega_K\left|B_K\right|\left|B_{j,K}\right|$
	for $j=1,\ldots,n$,
	where $B_{j,K}$ is obtained  from $B_K$ by replacing its $j$-th column with $(f_1, \dots, f_m)^\top$ for $j=1,\dots,n$.
	If $K\in\Pnws$, the interpolation conditions satisfy $B_K\bm{c}_K = \bm{f}_K$, where $\bm{f}_K=(f_{k_1}, \dots, f_{k_n})^\top$ with $K=\{k_1,\ldots,k_n\}$. Therefore, we write the solution to the linear problem as
	$\bm{c}_K=(c_{1,K}, \dots, c_{n,K})^\top$ which is unique for $\left|B_K\right|\neq 0$, and so by Cramer's rule we obtain
	$c_{j,K} = \left|B_{j,K}\right|\Big/\left|B_K\right|$.
	The associated interpolant $v_K(\bm x)$ for $K\in\Pnws$ can be written as $v_K(\bm x) = \sum_{j=1}^nc_{j,K}\beta_j(\bm x)$, hence by letting $\lambda_K=\omega_K\left|B_K\right|^2$, it holds
	\[
	\begin{split}
		v(\bm x) & = \sum_{j=1}^nc_j \beta_j(\bm x) = \sum_{j=1}^n\frac{\left|A_j\right|
			\beta_j(\bm x)}{\left|A\right|} = \sum_{j=1}^n\frac{\sum_{K\in \calP_{n}}\omega_K\left|B_K\right|\left|B_{j,K}\right|\beta_j(\bm x)}{ \sum_{K\in \calP_{n}}\omega_K\left|B_K\right|^2}\\
		& = \frac{\sum_{K\in \mathcal P^{\star}_{n}}\omega_K\left|B_K\right|^2\sum_{j=1}^n c_{j,K}\beta_j(\bm x)}{ \sum_{K\in \mathcal P^{\star}_{n}}\omega_K\left|B_K\right|^2} = \frac{\sum_{K\in \mathcal P^{\star}_{n}}\lambda_K  v_K(\bm x)}{ \sum_{K\in \mathcal P^{\star}_{n}}\lambda_K }.
	\end{split}
	\]
\end{proof}
\begin{remark}
	Note that \Cref{th:genweightedsum} is applicable to any finite-dimensional (multivariate) vector space. This encompasses various function spaces, such as the space of polynomials up to a specific degree, spline spaces with fixed degree and order, spline spaces with varying locations of knot lines or any other real function space of finite dimension equipped with a basis. Moreover, note that equation \eqref{eq:001} in \Cref{th:genweightedsum} is very similar to the formulation of multinode Shepard operators \cite{DELLACCIO2021} in the case of polynomial spaces $V$.
\end{remark}
\begin{remark}
	In general, we have {$\#\left(\calP_{n}^{\star}\right) \leq \#\left(\calP_{n}\right)=\binom{m}{n}$}, due to the non-nullity condition on the determinants $|B_K|$, which guarantees the existence of the interpolants. In particular, in the case of splines, the condition $|B_K|\neq0$ is equivalent to the Schoenberg-Whitney nesting condition \cite{deboor78}. In addition,
	the value of $\lambda_K$ depends on the location of the knot lines.
\end{remark}

\subsection{Consequences of the interpolatory formulation of weighted least squares approximation}\label{sec:consequences}
For the univariate polynomial case,  results on the upper and lower pointwise error bounds of the approximant derivatives up to a certain order as well as on the influence of the weights have been presented in \cite[Section 3]{floater2020} and \cite[Section 4]{floater2020}, respectively.
We extend these results to any vector space $V$ of finite dimension $n$. This generalization extends the derived conclusions of \cite{floater2020} in a broader setting, beyond the polynomial scenario.
More precisely, let $r\in \mathbb N$ be given and let  $\alpha = \left(\alpha_1, \ldots, \alpha_N\right)$ be any multi-index with $\sum_{j=1}^N\alpha_j = r \geq 0$. If the $r$-th order derivative of $v\in V$ exists, it can be expressed as the weighted average of the $r$-th order derivatives of the interpolants, i.e.,
	\begin{equation}\label{eq:der}
		\partial^\alpha v(\bm x) =  \frac{\sum_{K\in \mathcal P^{\star}_{n}}\lambda_K  	\partial^\alpha v_K(\bm x)}{ \sum_{K\in \mathcal P^{\star}_{n}}\lambda_K }.
	\end{equation}
In addition,  pointwise upper and lower bounds for the value of the $r$-th order derivative of $v(\bm x)$ can be obtained from \eqref{eq:der}, i.e.,
\begin{equation*}
	\min_{K\in \mathcal P^{\star}_{n}}\partial^\alpha v_K(\bm x) \leq \partial^\alpha v(\bm x) \leq \max_{K\in \mathcal P^{\star}_{n}}\partial^\alpha v_K(\bm x).
\end{equation*}
Likewise,  the pointwise approximation error shares the same weighted average, given by
\begin{equation*}
	f(\bm x)-v(\bm x) =  \frac{\sum_{K\in \mathcal P^{\star}_{n}}\lambda_K  \left(f(\bm x) - v_K(\bm x)\right)}{ \sum_{K\in \mathcal P^{\star}_{n}}\lambda_K }.
\end{equation*}
Furthermore, the following consequences on the influence of the weights $\omega_i$ for $i=1, \ldots, m$ and derivative estimations can be inferred.
\begin{corollary}\label{cor:weight}
	Let  $I\subseteq \{1, \ldots, m\}$ be of cardinality $\#(I) = r$, with $1 \leq r \leq n$. Define
	\[
	u_I(\bm x) = \lim_{\substack{{\omega_i \to +\infty}\\ \forall i\in I}} v(\bm x)\quad \text{and}\quad \calP_{n}^{\star}\setminus I = \{ K\subseteq \{1,\dots, m\}\setminus I\colon \#(K)=n,\ \text{and}\ \left\vert B_K\right\vert \neq 0 \}.
	\]
	If $r=n$, then $u_I(\bm x) = v_I(\bm x)$, represents the interpolant of the data points indexed by $I$. 
	If $1\leq r < n$, then  
	\[
	u_I(\bm x) = \frac{\sum_{K\in \mathcal P^{\star}_{n-r}\setminus I}\lambda_{I,K}  v_{I\cup K}(\bm x)}{ \sum_{K\in \mathcal P^{\star}_{n-r}\setminus I}\lambda_{I,K}},\quad \text{with}\quad \lambda_{I,K} = \omega_K \left\vert B_{I\cup K}\right\vert^2.
	\]
\end{corollary}
\begin{proof}
	The proof can be obtained by  by considering $V = \spn\{\Gamma\}$, where $\Gamma = \{\beta_1, \ldots, \beta_n\}$ with $\beta_j : \Omega \to \mathbb{R}$ for $j=1, \ldots, n$, and  substituting the Vandermonde matrices in \cite{floater2020} with the corresponding collocation matrices $B$ in \eqref{eq:colloc} and $B_K$ in \eqref{eq:colloc-loc}.
\end{proof}

Another important implication of~\Cref{th:genweightedsum} is the possibility of rewriting $\ell^p$--approximation problems as suitable convex combination of interpolants as outlined in \Cref{rem:irls} where we exploit the Iterative Reweighted Least Squares (IRLS) method, see e.g.,~\cite[Section~4.5]{Bjorck1996}. 
\begin{remark}\label{rem:irls}
Let us consider the problem
\begin{equation}\label{eq:lp_problem}
\min_{u\in V} \sum_{i=1}^m { \left\| u(\bm x_i) - f_i\right\|_p^p}
\end{equation}
with $1<p<2$, {where} $V = \spn\left\{\Gamma\right\}$ {is} a vector space of dimension $n$, and $\Gamma = \left\{\beta_1,\ldots,\beta_n\right\}$ its basis. 
The IRLS method approximates the exact solution of problem~\eqref{eq:lp_problem} through the iterative computation of weighted least squares problems. 
In particular, the weights are recursively updated for a maximum number of iterations $K_{\mathrm{max}}$ {by}
\[
\omega_i^{k+1} = \left\vert u^k(\bm{x}_i) - {f}_i\right\vert^{\nicefrac{(p-2)}{2}},
\]    
with weights $\omega_i^{k}$ and solution  $u^k$ of problem~\eqref{eq:Vwls}, for {iterations} $k=1, \ldots, K_{\mathrm{max}}$. 
By direct application of~\eqref{eq:001} in~\Cref{th:genweightedsum}, the outcome of each iteration can be rewritten as a convex combination of interpolants.
Note that the interpolants and the determinant of the matrices $B_K$ do not need to be recomputed. 
In other words, it is enough to update the weights $w_K$ {according to} \eqref{eq:prodwk}, to compute the solution of problem~\eqref{eq:lp_problem}.
\end{remark}

\section{Hierarchical spline models}
\label{sec:splines}
The  specification of the basis used to design  
the continuous model $\bm v\colon \Omega \to \mathbb{R}^D$ approximating the observations $\{\left(\bm x_i, \bm f_i\right)\}_{i=1}^m$, with $\bm x_i\in\Omega\subset\mathbb{R}^N$ and $\bm f_i\in\mathbb{R}^D$ for $i=1, \ldots, m$,  plays a fundamental role both for the geometrical and numerical properties of the approximant $\bm v$. Beyond polynomial constructions, splines have been largely adopted both for geometric modelling \cite{farin1993} and isogeometric analysis \cite{hughes2005} because of their properties as well as the simplicity of their constructions.
In particular, for $N=1$, let $\Omega = [a, b]\subset \mathbb{R}$, and $\left\{a \equiv \tau_0 < \tau_1 < \ldots < \tau_{L-2} < \tau_{L-1} \equiv b\right\}$ a partition of $\Omega$, which defines $ \gamma_j = [\tau_{j}, \tau_{j+1})$ for $j=0, \ldots, L-3$ and $\gamma_{L-2} = [\tau_{L-2}, \tau_{L-1}]$. Moreover, set a certain polynomial degree $d$ and corresponding order  $k = d+1$ and
the multiplicities  $\left\{\mu_0, \ldots, \mu_{L-1}\right\}$ associated to each $\tau_j$, with $1 \leq \mu_j\leq k$, for $j=0, \ldots, L-1$. These elements lead to the definition of the so called \emph{knot vector} $\bm t = [\tau_0, \tau_1, \ldots, \tau_{L-1}]$, whose items $\tau_j$ appears according to their multiplicity $\mu_j$, for $j=0, \ldots, L-1$, and finally to the polynomial  spline space on $\Omega$, i.e.
\[
\left\{ u: \Omega \to \mathbb{R} \quad \vert \quad  u_{|\gamma_j}\in \Pi_{k},\; u^{(r)}(\tau_j^{-}) =  u^{(r)}(\tau_j^{+})\;  \text{for}\; r=1,\ldots,k-\mu_j \; \text{and}\; j=0, \ldots, L-2 \right\},
\]
where $\Pi_{k}$ indicates the space of polynomials of maximum order $k$ and $u^{(r)}$ indicates the $r$-th order derivative of $u$.  Moreover, univariate splines are non-negative, have local support and form a partition of unity.

When moving to multivariate settings, for $N>1$, hierarchical splines provide a natural strategy to preserve locality and therefore develop adaptive spline constructions with local refinement capabilities. In the context of data fitting, {this} allows to
suitably design adaptive spline approximation schemes.
For $N>1$, we introduce the  polynomial multi-order $\bm k$ and consider a sequence of \emph{nested} linear vector \emph{spline} spaces defined by
$V^0 \subset \ldots \subset V^\ell \subset V^{\ell+1} \subset \ldots \subset V^{L-1}$, where the index $\ell$ of $V^\ell$ is called its level.
{We assume that} {$V^\ell~=~\spn\left\{\Gamma^\ell_{\bm k}\right\}$ have finite dimension $n_\ell = \mathrm{dim}V^\ell$ {where} the polynomial spline basis $\Gamma^\ell_{\bm k} = \left\{\beta_1^\ell, \ldots, \beta_{n_\ell}^\ell\right\}$ {is} defined by the polynomial multi-order $\bm k$ and spline basis functions $\beta_j^\ell:\Omega \to \mathbb{R}$. {In addition,} for each level $\ell$, we assume that the boundaries of the supports of all basis functions $\beta_j^\ell$, for $j=1, \ldots, n_\ell$, partition the
domain $\Omega$ into a {certain} number of connected \emph{cells} of level $\ell$. The collection of all cells among all the levels defines a tessellation $T(\Omega)$ of the domain $\Omega$. 
In addition, let $\Omega \equiv \Omega^0 \supset \ldots \supset \Omega^\ell \supset \Omega^{\ell+1} \supset \Omega^{L}\equiv \emptyset$ be a sequence of nested open subdomains, where each $\Omega^\ell$ represents a local region of $\Omega$ where additional degrees of freedom are required. In particular, we assume that the closure of each subdomain $\Omega^\ell$ coincides with the closure of a collection of cells of level $\ell$. 
The hierarchical spline construction consists of replacing any spline basis function of level $\ell$ whose support is completely contained in $\Omega^{\ell+1}$ by splines at successively hierarchical levels. More precisely, the hierarchical spline basis is defined as
\begin{equation}\label{eq:hspline}
\mathcal{H}_{\bm k} \coloneqq
\left\{ \beta_j^\ell \; \vert \; j\in A^\ell_{\bm k},\, \ell = 0, \ldots, L-1 \right\}
\end{equation}
with the indices
\[
A^\ell_{\bm k} \coloneqq
\left\{j\in\Gamma^\ell_{\bm k} \; \vert \; \text{supp}\left(\beta_j^\ell\right) \subseteq \Omega^\ell \land \text{supp}\left(\beta_j^\ell\right) \not\subseteq \Omega^{\ell+1} \right\}
\] where $\text{supp}\left(\beta_j^\ell\right)$ denotes the intersection of the support of $\beta_j^\ell$ with $\Omega^0$. The corresponding hierarchical space is defined as $\spn\left\{\mathcal{H}_{\bm k}\right\}$.

Several families of spline basis functions and their corresponding function spaces admit a hierarchical extension. 
That is the case of B-splines \cite{deboor78}, NURBS \cite{piegl1996}, as well as box splines \cite{deboor1993}, among others.
Examples of their hierarchical formulations can be found in \cite{forsey1988,kraft1997,kang2015}. 
In the following, we provide examples for \Cref{th:genweightedsum} and \Cref{cor:weight} in the special case of B-spline and hierarchical box-spline models, respectively.

\subsection{Verification of \Cref{th:genweightedsum} for spline models}
In this section, we present {the} numerical verification of  \Cref{th:genweightedsum} for two different vector spaces $V$: the space of polynomials and the space of univariate splines with a specific degree and regularity.
We consider a set of $m=7$ observations in the form ${(x_i,f_i)}$, for $i=1, \ldots, m$, given by
$\{(-4.5, -2), (-3.5, 0), (-2.2, -1), (-1.2, 2.8), (0.8, 2.9), (2.2, 0.5), (4.0, -2)\}$.
Furthermore, we assign weights $\omega_i$ for $i=1, \ldots, m$ {uniformly distributed on $(0,1)$}.
In the polynomial case, we choose the polynomial degree $d = 2$,
which implies ${\#\left(\calP_{d+1}\right)} = \binom{m}{d+1} = \binom{7}{3} = 35$. Thus, we have a total of 35 interpolation problems to be solved. Moving on to polynomial spline spaces, in addition to the degree, we {consider the spline order $k=d+1=3$ and set}
the knot vector $\bm{t} = [-5, -5, -5, -5/3, 5/3, 5, 5, 5]$ of length 8, {implying that the spline space has dimension $n = 8 - 3 =  5\leq m$}.
{Since $n\leq m$, there are} at most ${\#\left(\calP_{n}\right)} = \binom{m}{n} = \binom{7}{5} = 21$ interpolation problems required to fully reconstruct the spline least squares approximant.
In the case of spline spaces, there is no guarantee that each point subsequence satisfies the Schoenberg-Whitney nesting conditions. To illustrate this, consider the data set identified by $K=\{1,2,3,4,5\}\subset \{1,\ldots,m\}$, which does not satisfy these conditions. This is due to the absence of {data points} within the support of the last basis function, specifically $x_i\not\in[5/3, 5]$ for $i=1,\ldots, 5$. Consequently, the interpolant $v_K$ does not exist, and no interpolation problem needs to be solved for this particular subset.
To fully reconstruct the final global least {squares} approximation in this example, we need to compute 20 interpolation problems instead of the original 21.
\Cref{fig:sK5} illustrates the interpolants and the global least squares approximation involved in \eqref{eq:001} for the polynomial and spline spaces introduced for this numerical example. Note that we choose to show an example related to the univariate polynomial and spline space, but our results hold in general for any (multivariate) vector space of finite dimension.

\begin{center}
	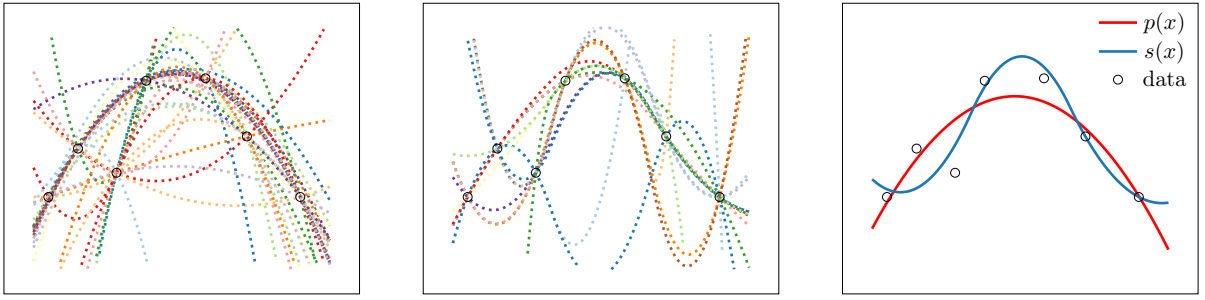
\begin{figure}
		\centering
		\resizebox{\linewidth}{!}{
			\begin{tikzpicture}
				\begin{groupplot}[
					group style={
						group name=group,
						group size=3 by 1,
						xlabels at=edge bottom,
						xticklabels at=edge bottom,
						ylabels at=edge left,
						yticklabels at=edge left,
					},
					ymin = -5,
					ymax = 5,
					ymajorticks=false,
					xmajorticks=false,
					restrict y to domain=-5:5,
					width=0.45\linewidth,
					enlargelimits
					]
					
					\pgfplotsset{cycle list/Paired-12}
					\nextgroupplot[cycle multi list={Paired-12}]
					\foreach \x in {1,...,35}
					{
						\addplot+[dotted, very thick] table [x=x,y=p\x, col sep=comma]{polynomial.txt};
					}
					\addplot+[mark = o, only marks, mark options={draw=black, fill = black}] coordinates{
						(-4.5000,   -2.0000)
						(-3.5000 ,        0)
						(-2.2000  , -1.0000)
						(-1.2000   , 2.8000)
						(0.8000 ,   2.9000)
						(2.2000 ,   0.5000)
						(4.0000 ,  -2.0000)
					};
					\pgfplotsset{cycle list/Paired-12}
					\nextgroupplot
					\foreach \x in {1,...,16}
					{
						\addplot+[dotted, very thick] table [x=x,y=s\x, col sep=comma]{splines.txt};
					}
					\addplot+[mark = o, only marks, mark options={draw=black, fill = black}] coordinates{
						(-4.5000,   -2.0000)
						(-3.5000 ,        0)
						(-2.2000  , -1.0000)
						(-1.2000   , 2.8000)
						(0.8000 ,   2.9000)
						(2.2000 ,   0.5000)
						(4.0000 ,  -2.0000)
					};
					\pgfplotsset{cycle list/Set1}
					\nextgroupplot[legend style={fill=none, draw=none, at={(1,1)}, anchor=north east}]
					\addplot+[very thick, draw=red] table [x=x,y=plocal, col sep=comma]{solutions_polynomial.txt};\addlegendentry{$p(x)$}
					\addplot+[very thick] table [x=x,y=slocal, col sep=comma]{solutions_spline.txt};\addlegendentry{$s(x)$}
					\addplot+[mark = o, only marks, mark options={draw=black}] coordinates{
						(-4.5000,   -2.0000)
						(-3.5000 ,        0)
						(-2.2000  , -1.0000)
						(-1.2000   , 2.8000)
						(0.8000 ,   2.9000)
						(2.2000 ,   0.5000)
						(4.0000 ,  -2.0000)
					};
					\addlegendentry{data}
				\end{groupplot}
			\end{tikzpicture}
		}
		\caption{Numerical verification of \eqref{eq:001}, \Cref{th:genweightedsum}: $v_K(\bm x)$ interpolants for polynomials (left) and splines (center) and the final weighted least {squares} approximations as sum of interpolants (right) for polynomials (red) and splines (blue).}
		\label{fig:sK5}
	\end{figure}
\end{center}

\subsection{Weighted least squares with hierarchical box splines}\label{exm:box-splines}
In this section, we perform a numerical investigation of  the role of the weights addressed in  \Cref{cor:weight} in \Cref{sec:consequences} 
for a weighted least squares problem with hierarchical box splines in the bivariate case.
We consider a point cloud of $64\times 64$, i.e. $m = 4096$, uniformly gridded data, obtained by sampling the function
\begin{equation}\label{eq:3peaks}
\begin{split}
	& f(\bm x) = f(x, y) = \\
	& \frac{2}{3 \exp \left(\sqrt{\left(10x-3\right)^2 + \left(10y-3\right)^2}\right)} + \frac{2}{3\exp \left(\sqrt{\left(10x+3\right)^2 + \left(10y+3\right)^2}\right)} + \frac{2}{3\exp\left(\sqrt{\left(10x\right)^2 + \left(10y\right)^2}\right)}
\end{split}
\end{equation}
for $\bm x = (x_1, x_2)\in[-1, 1]^2$ and we associate to each item of the point cloud a weight $\omega_i > 0$ for $i=1, \ldots, m$. 
For different weight values, we compute the weighted least {squares} approximation $v\left(\bm x\right)$ of the input data in terms of 521 $C^2$ quartic box splines, in their hierarchical extension, see again \cite{kang2015}. In order to investigate the role of the weights in the approximation process, we execute the following two steps.
\begin{enumerate}[(a)]
	\item Firstly, we perform a standard least {squares} approximation and individuate the data sites whose approximation error is above a certain threshold $\epsilon$, thus define 
	$$K \coloneqq \{ i \in \{1, \ldots, m\} : \left\vert v(\bm x_i) - f(\bm x_i)\right\vert >~\epsilon\}.$$  
	\item Subsequently, for $i\in \{1, \ldots, m\}\setminus K$, we set the corresponding weights values $\omega_i = \omega_0$ and similarly for $i\in K$ we set $\omega_i = \omega_\gamma$.
\end{enumerate}
Therefore, we analyse the accuracy behaviour of the resulting approximant $v(\bm x)$ for different choices of $(\omega_0, \omega_\gamma)$ in terms of maximum approximation error (MAX). 
Setting $\epsilon = 5$e$-4$, then $\#(K) = 124$ and for $\omega_0 = \omega_\gamma = 1$, the resulting ordinary least {squares} approximation is characterized by a MAX error of 2.15e$-$3.
We then keep $\omega_0 = 1$ and vary $\omega_\gamma = 2, 3, 4, \ldots, 10, 50, 100$.
We note that increasing the value of the weights $\omega_\gamma$ may help the final accuracy of the approximant. 
In particular, we obtain MAX = 2.15e$-$3, 2.06e$-$3, 2.00e$-$3 for $\omega_\gamma= 2, 3, 4$ and MAX decreases  further when increasing $\omega_\gamma$ until achieving its minimum, with MAX =  1.89e$-$3, for $\omega_\gamma = 6$. 
However, if the weight values $\omega_\gamma$ are too large, the final approximant deteriorates. 
For instance, for $\omega_\gamma = 7$, MAX = 1.91e$-$3, for  $\omega_\gamma = 10$, MAX = 2.63e$-$3 and for  $\omega_\gamma = 100$, MAX = 1.36$-$2. 
These results can also be observed in \Cref{fig:c2hboxapprox}, which depicts the hierarchical box spline mesh (left) together with two weighted least {squares} approximation resulting from two different choices of weights, namely $\omega_\gamma = 6$ (center) and $\omega_\gamma = 100$ (right). 
The pointwise error color map, ranging from the minimum to the maximum approximation error,  is plotted for each surface separately. The corresponding color bars are also reported.

As shown in this example, the choice of weights can lead to either better or worse spline fitting results in terms of the MAX error compared to ordinary least squares, where all the weights are set to ones. 
To the best of our knowledge, only a few deterministic methods address the use of weights associated with data points in the spline fitting literature. However, the topic is widely discussed in the statistics community, going back to locally weighted scatterplot smoothing (LOWESS)~\cite{Cleveland1981} for smoothing scatterplots by robust locally weighted regression \cite{Cleveland1979} and multivariate adaptive regression splines (MARS) \cite{Friedman1991}, with a more recent method proposed in \cite{BRUGNANO2024} and the references therein, for example. In this work, we propose a deterministic strategy to suitably take advantage of using weights associated with data points within the spline fitting process which also includes adaptive spline fitting constructions.

\begin{figure}[!t]
	\centering
	\resizebox{\linewidth}{!}
	{
		\begin{tikzpicture}
			\centering
			\node (a) {\includegraphics[width=0.3\linewidth, trim=130 250 110 220, clip]{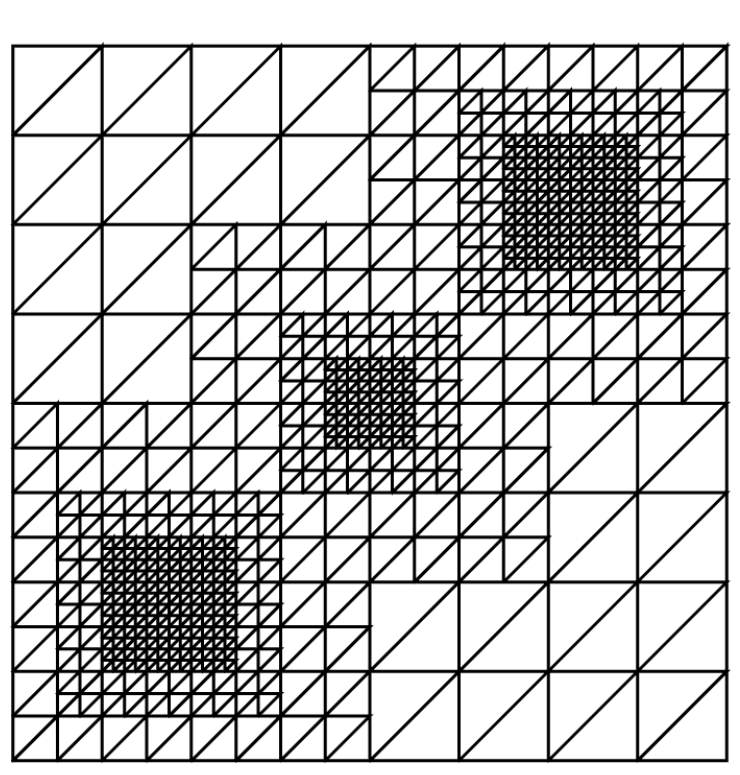}};
			\node (b) [right = 0pt of a] {\includegraphics[width=0.3\linewidth, trim = 160 250 80 220, clip]{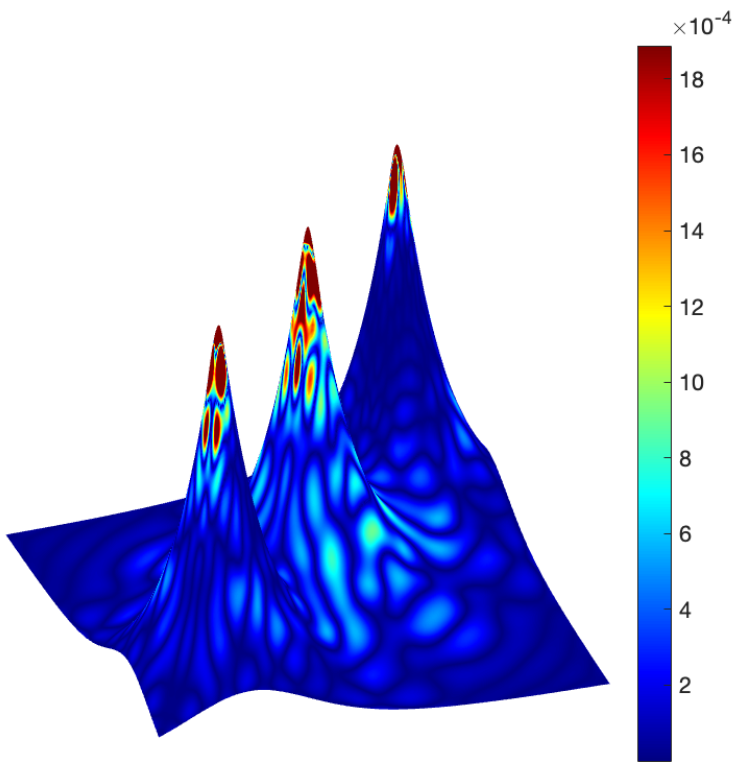}};
			\node (c) [right = 0pt of b] {\includegraphics[width=0.3\linewidth, trim = 160 250 80 220, clip]{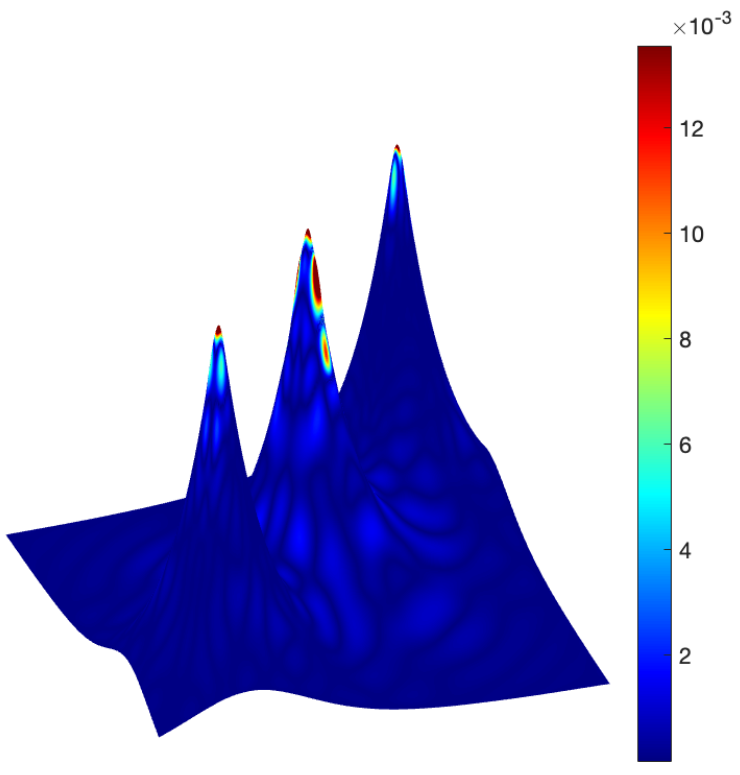}};
			\node (bt) [below =0pt of b] {\footnotesize $\omega_\gamma =  6$, MAX = 1.89e$-$3};
			\node (ct) [below =0pt of c] {\footnotesize $ \omega_\gamma = 100$, MAX = 1.36e$-$2};
		\end{tikzpicture}
	}
	\caption{Hierarchical box spline mesh (left); weighted least {squares} surfaces for $\omega_\gamma=6$ (center) and 100 (right).}\label{fig:c2hboxapprox}
\end{figure}

\section{Reweighted least squares spline fitting}
\label{section:rWLS}
In this section, we propose a strategy to take effective advantage of the weights associated with {each} input data in the context of fitting problems. 
This approach is inspired by the notion of landmarks used in shape analysis~\cite{Bookstein1978}, but
we introduce a more general concept of \emph{markers} and use them within the framework of curve and surface fitting. 
More precisely, landmarks can be understood as a set of labelled points which represent some physically identifiable parts of an object,
as well as important features of the input data, which need to be encoded and reproduced in the final continuous approximation model. 
\Cref{fig:landmarks} shows four point clouds (in blue) and the chosen landmarks (in black), which 
represent the features desirable to be preserved by the fitted curve. 
However, depending on the acquisition process, data can be affected by noise and outliers whilst the final approximation models should avoid the reproduction of corrupted data.  
For  the given set of observations  $\{\bm x_i,\bm f_i\}_{i=1}^m$, we generalize the concept of landmarks to \emph{markers} of  two {types}. More precisely, we define the index set $K_I \subseteq\{1,\ldots,m\}$ as \emph{markers of type~I} if the associated points $\{\bm x_i,\bm f_i\}$  for $i\in K_I$  represent data features {to be preserved}, while we define the index set $K_{II} \subseteq\{1,\ldots,m\}$ as \emph{markers of type~II} if the index indicates noisy data or outliers which should not be reproduced. 
In particular, we have $K_I \cap K_{II} = \emptyset$ and $K_I \cup K_{II} \subseteq \{1, \ldots, m\}$. 
Note that the choice of type I and type II markers and their identification depends considerably on the problem at hand and it is still an open research topic, see e.g., for automatic features selection \cite{zulqarnain2015shape, gupta2015knowledge} or for outliers recognition \cite{ning2018efficient,farahmand2012robust,farahmand2011doubly}. 
The markers identification is out of the scope of the present article and we assume the markers (of both types) to be known a priori. Nevertheless, for synthetic data, we devise an error-driven detection of markers second to point cloud pre-processing.
Finally, we provide fitting schemes which can address simultaneously both types of markers, by leveraging the weight values associated with them depending on the approximation error, also within an adaptive approximation framework.

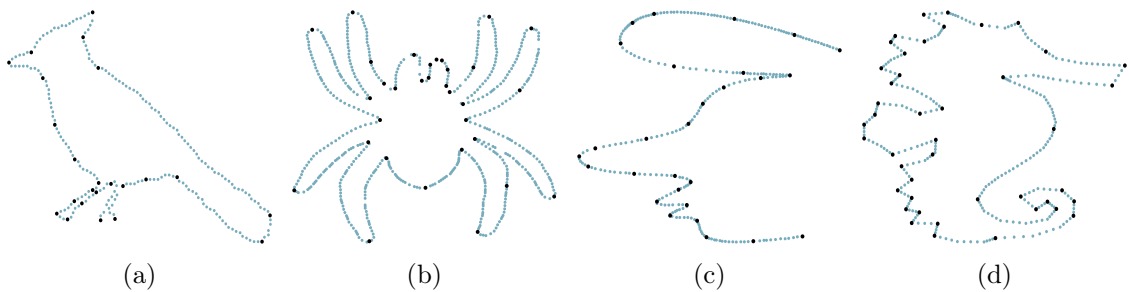
\begin{figure}[!t]
	\begin{center}
		\begin{tikzpicture}
			\node (bird) {\includegraphics[width=0.22\textwidth,height=0.15\textheight]{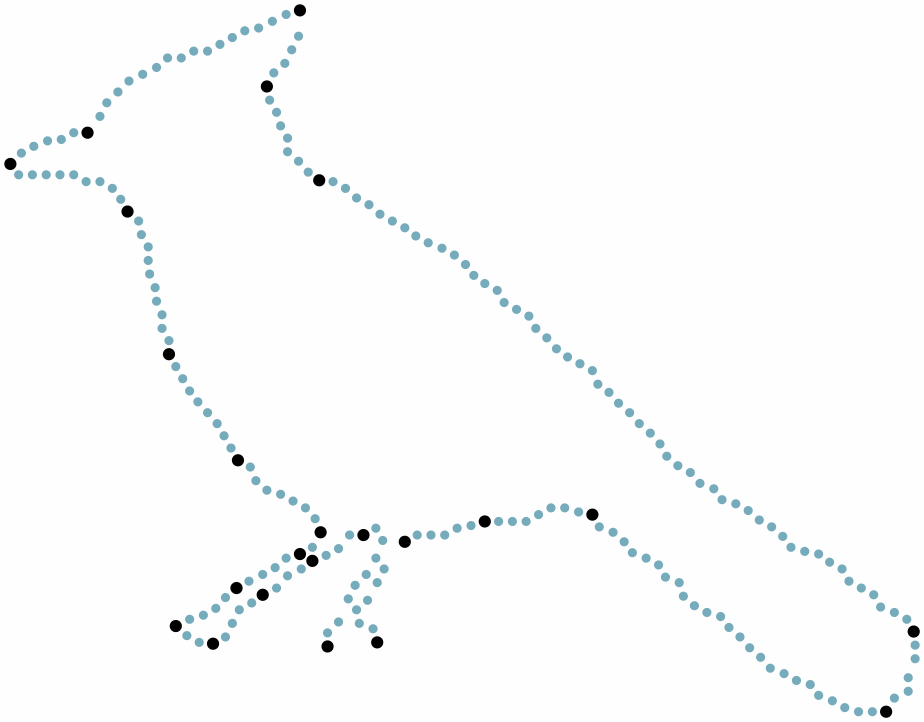}};
			\node (spider) [right = 0pt of bird] {\includegraphics[width=0.22\textwidth,height=0.15\textheight]{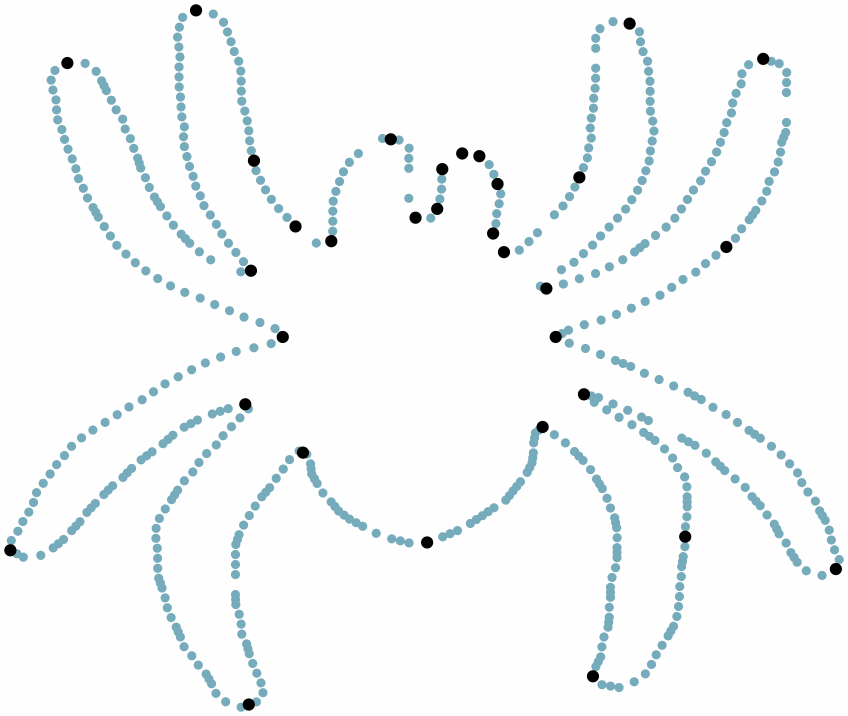}};
			\node (face) [right = 0pt of spider] {\includegraphics[width=0.22\textwidth,height=0.15\textheight]{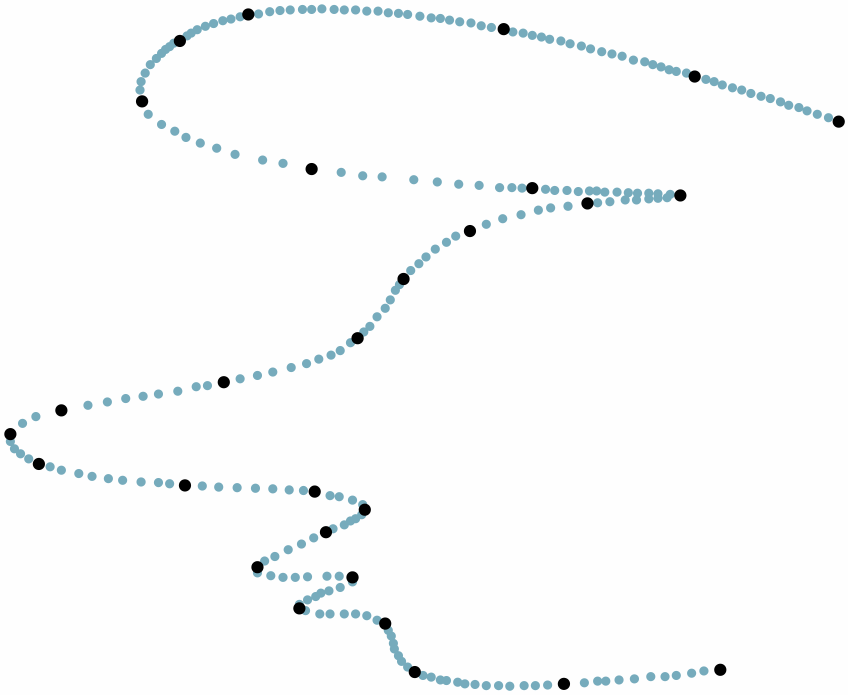}};
			\node (seahorse) [right = 0pt of face] {\includegraphics[width=0.22\textwidth,height=0.15\textheight]{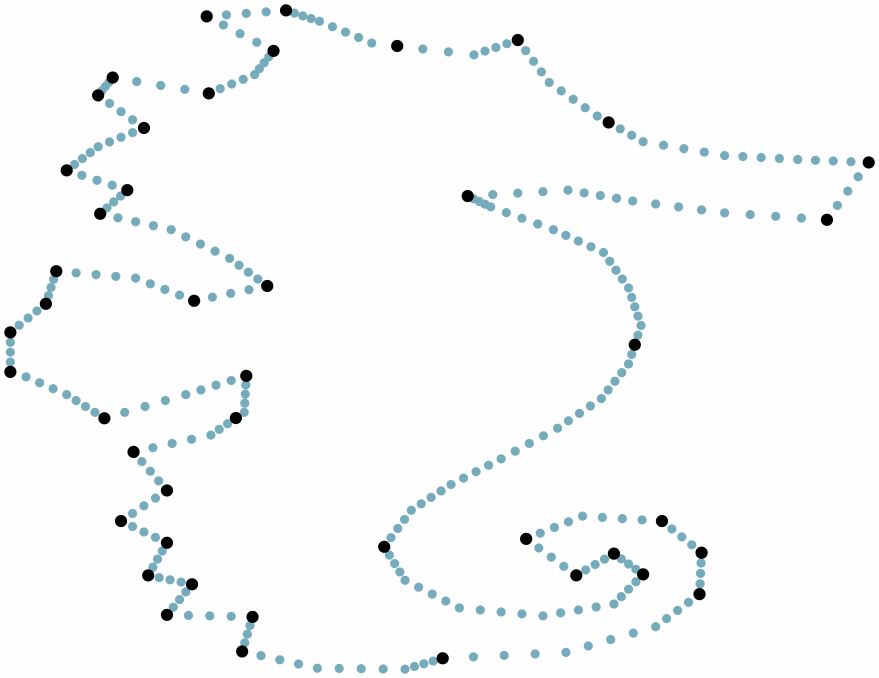}};
			\node [below = 0pt of bird] {(a)};
			\node [below = 0pt of spider] {(b)};
			\node [below = 0pt of face] {(c)};
			\node [below = 0pt of seahorse] {(d)};
		\end{tikzpicture}	
		\end{center}
	\caption{Point cloud (blue dots) and set of type I markers (black dots)}
	\label{fig:landmarks}
\end{figure}

The fitting problem with markers consist of finding an element $\bm v\in V$ which approximate the points  $\{\bm x_i,\bm f_i\}_{i=1}^m$, i.e. $\bm v(\bm x_i) \approx \bm f_i$, with
$\|\bm v(\bm x_i) - \bm f_i\|_2 < \mathrm{tol}_I$ for $i\in K_I$ and $\|\bm v(\bm x_i)  - \bm f_i\|_2 < \mathrm{tol}_{II}$ for $i\not\in K_{II}$.
Note that the usual formulation of a fitting problem can be interpreted as a special instance of the present one. 
In particular, it is equivalent to setting $K_I = \left\{1, \ldots, m\right\}$, $K_{II} = \emptyset$ and choosing $\mathrm{tol}_I = \mathrm{tol}_{II} = \epsilon$.
Specifically, the reweighted least squares algorithm with markers is described in~\Cref{alg:rWLS-landmarks} and it consists of the following steps. 
For an initial choice of the weight values, 
\begin{enumerate}
	\item\label{step:rWLS_1} 
	Solve the weighted least squares approximation problem \eqref{eq:Vwls}.
	\item\label{step:rWLS_2} 
	Update the point-wise error $e_i = \|\bm v(\bm x_i) - \bm f_i\|_2$ for each $i=1, \ldots, m$.
	\item\label{step:rWLS_3}
	Check whether the current fitting $\bm v\in V$ meets the requirements prescribed by the markers $K_I$ and $K_{II}$ and the respective error tolerances. In particular, the fitting process is completed if $e_i < \mathrm{tol}_I$ for $i \in K_I$ and $e_i < \mathrm{tol}_{II}$ for $i\not\in K_{II}$.
	\item\label{step:rWLS_4}
	If the accuracy requirements are satisfied, return the approximation computed in \Cref{step:rWLS_1},~otherwise update the weight values corresponding to the markers. In particular, $\omega_i = \omega_i \cdot \alpha$ with $\alpha > 1$ for each $i\in K_I$ and $\alpha < 1$ for each $i\in K_{II}$.
	\item\label{step:rWLS_5} 
	Start again from \Cref{step:rWLS_1} with the new weight values.
\end{enumerate}
Note that the update choice of the weights in \Cref{step:rWLS_4} increases the values of the weights related to markers of type I, whereas it decreases them for markers of type II. 
In particular, we suggest an error-driven update of the weights with $\alpha = (1 + e_i)$ for $i\in K_I$ and $\alpha = {1}/{(1 + e_i)}$ for $i\in K_{II}$. 
Finally, it is worth highlighting that the IRLS procedure \cite{osborne1985,Bjorck1996} falls within this general fitting scheme for the special choice of $\alpha = {1}/{\max\{\delta, e_i\}}$ for each $i\in K_{II}$, with $\delta>0$ necessary for the stability of the IRLS method.
\begin{algorithm2e}[!t]
	\SetAlgoLined
	\KwData{Point cloud $\{\bm x_i,\bm f_i\}_{i=1}^m$, the set of markers $K_I, K_{II}\subseteq\{1,\ldots,m\}$, the tolerances $\textup{tol}_I$, $\textup{tol}_{II}$, a fixed vector space $V = \spn\{\beta_1, \ldots, \beta_n\}$ and a maximum number of iterations $\mathrm{M}_\mathrm{max}$;}
	\KwResult{$\bm{v} \in V$ the reweighted least squares approximant.}
	Initialize the weights $\omega_i = 1$ and the point-wise errors $e_i = 1$ for each $i=1, \ldots, m$ and $\mathrm{loop}=0$\\
	\While{$\max_{i\in K_I}{e_i} > \textup{tol}_I$ and $\max_{i\in I\setminus K_{II}}{e_i} > \textup{tol}_{II}$ and $\mathrm{loop}<\mathrm{M}_{\mathrm{max}}$}{%
		Solve the weighted least squares problem \eqref{eq:Vwls};
		\\
		compute the errors $e_i = \|\bm{v}(\bm x_i) - \bm f_i \|_2 	$; \\
		update the weights associated to the landmaks $K_I$ and $K_{II}$, as $w_i = w_i \cdot\alpha(e_i)$;
		\\
		set loop = loop$+$1;
	}
	\KwRet{$\bm{v} \in V$ which  solves the weighted least squares problem.}
	\caption{General formulation of the reweighted least squares fitting.}
	\label{alg:rWLS-landmarks}
\end{algorithm2e}
The reweighted fitting scheme from \Cref{alg:rWLS-landmarks} is presented in its more general formulation, namely, the approximant is sought in any \emph{fixed} finite-dimensional vector space $V$.
However, in the incoming numerical experiments, we apply it to spline models as the ones defined in \Cref{sec:splines}.

As far as spline constructions are concerned, the definition of a fixed spline space $V = \spn\left\{\beta_1, \ldots, \beta_n\right\}$ relies on the choice of the polynomial degree and order as well as a suitable tessellation $T(\Omega)$ of the domain $\Omega$. 
For spline fitting problems, in addition to the definition of the spline space $V$ also the parameterization of the data needs to be addressed, namely to associate at each $\bm f_i\in\mathbb{R}^D$ a suitable parametric value $\bm x_i\in \Omega$.
Both \emph{knot placement} to define $T(\Omega)$ and \emph{parameterization} to define the sites $\bm x_i\in \Omega$ for $i=1, \ldots, m$, play fundamental roles in the accuracy of the final spline fitting model, they are crucial open research topic and several methods have been provided to properly address them, see e.g. \cite{lee1989,floater2001,wang2013,zhang2016,fang2013,hu2020,zhu2022} and the references therein. 
Their resolution is out of the scope of the present paper, thus we assume both the parameterization and the spline space to be constructed with established techniques. 
Nevertheless, special attention needs to be dedicated to adaptive spline construction, where the space $V$ is not fixed but enlarged and updated iteratively.

\subsection{Reweighted least squares adaptive spline fitting}\label{subsec:adaptive-rWLS}
In this section, we show how to suitably combine the update of the weights with adaptive spline approximation schemes.
We revisit the adaptive least squares fitting scheme proposed in \cite{giannelli2023eg} by suitably assigning weight values to the data observation within the adaptive routine. 
The main idea which drives an adaptive fitting algorithm consists of adding iteratively degrees of freedom in regions of the domain $\Omega$ where the approximation error is too high. 
Similarly to \Cref{alg:rWLS-landmarks}, if the points with a too high error belong to $K_I$, then their weights will be augmented, otherwise if they belong to $K_{II}$ their weights will be diminished. 
In addition, at each iteration of the adaptive loop, not only the updates of the weight values take place, but also of the sets $K_I$ and $K_{II}$. 
This strategy is effective since, thanks to the adaptive refinement, in some regions of the domain the accuracy requirements $\mathrm{tol}_I, \mathrm{tol}_{II}$ are already locally achieved and it is useless or even harmful to keep on modifying the weight values.
More precisely, once an initial parameterization and tessellation $\Omega$ are chosen, i.e. for a fixed hierarchical spline space $V$, and for an initial choice of the weights we perform the following steps.
\begin{enumerate}
\item\label{step:rWLSAdap_1} 
Compute the (penalized) weighted least squares problem
	\begin{equation}\label{eq:wls-pen}
		\bm{v}(x) = \argmin_{\bm u \in V}\frac{1}{2}\sum_{i=1}^{m}\omega_i\left\| \bm u\left(\bm x_i\right) - \bm f_i\right\|_2^2 + \lambda J\left(\bm u\right),
	\end{equation}
	where the penalization term $J$ is the thin-plate energy functional, whose influence is controlled by $\lambda \geq 0$, i.e. for $\bm x = (s,t)\in\Omega$,
	\begin{equation}\label{eq:pen-fun}
	J\left(\bm u\right) = \bigintsss_{\Omega}\left\|\frac{\partial^2\bm u}{\partial s\partial s}\right\|_2^2 + 2\left\|\frac{\partial^2\bm u}{\partial s\partial t}\right\|_2^2 + \left\|\frac{\partial^2\bm u}{\partial t\partial t}\right\|_2^2 \mbox{d}s\mbox{d}t.
	\end{equation}
\item \label{step:rWLSAdap_2} 
Evaluate the error indicator as the point-wise error distance $\left\|\bm v(\bm x_i) - \bm f_i\right\|_2$ for each $i=1, \ldots, m$. By individuating the sites {$\bm x_i \in \Omega$} where $\left\|\bm v(\bm x_i) - \bm f_i\right\|_2 \geq \epsilon$, the error indicator identifies the region of the domain $\Omega$ where potentially additional degrees of freedom are needed to meet the prescribed accuracy $\epsilon$.

\item \label{step:rWLSAdap_3} 
If $\left\|\bm v(\bm x_i) - \bm f_i\right\|_2 < \epsilon$ for each $i=1, \ldots, m$ return the approximation computed in \Cref{step:rWLSAdap_1}, otherwise update the weights valued and the markers $K_I$ and $K_{II}$. 
In particular, for each $i\in K_I$ if $e_i < \mathrm{tol}_I$, then the $i$-th data satisfies the accuracy requirements and it should not belong to the markers of type I any more, namely $K_I = K_I\setminus\{i\}$. 
On the contrary, if $e_i > \mathrm{tol}_I$, then its corresponding weight needs to be enlarged, i.e. $\omega_i = \alpha\cdot\omega_i$, with $\alpha > 1$. 
Similar considerations hold for $i\in K_{II}$, with inverted inequalities and $\alpha < 1$.

\item \label{step:rWLSAdap_4} 
Mark the cells of the current hierarchical spline space which contain the parameters $\bm x_i$ identified by the error indicator.

\item \label{step:rWLSAdap_5} 
Refine the marked cells by suitably splitting them to effectively enlarge the hierarchical spline space.

\item \label{step:rWLSAdap_6} 
Update the data parameterization by suitably moving the data sited $\bm x_i$ for each $i=1, \ldots,m$.

\item \label{step:rWLSAdap_7} 
Start again from \Cref{step:rWLSAdap_1} with the new weight values.
\end{enumerate}

Note that in  \Cref{step:rWLSAdap_1} we introduced the penalization term \eqref{eq:pen-fun}, usually addressed as \emph{thin-plate energy}, whose influence is ruled by a weight $\lambda\geq0$.
The introduction of such a functional is a common practice in spline geometric modelling, 
in particular when reconstructing a spline geometry from a set of unorganized data.
More precisely, a regularization term is commonly introduced to smoothen the solution and therefore avoid the presence of spurious oscillations and artefacts, which may affect the final geometric model.
According to \cite{kiss2014}, we set $\lambda\geq0$ to be a constant small ($\sim10^{-6}$) value and we keep it fixed during the entire fitting process. 
Non-constant regularization weight functions for data fitting via least-squares tensor-product splines have been recently proposed in \cite{merchel2023,lenz2023}, see also the references therein, e.g., \cite{craven1978, wahba1990, gu1992}.

For more details on \Cref{step:rWLSAdap_4,step:rWLSAdap_5,step:rWLSAdap_6}, about marking refinement and parameterization routines, we refer the reader to \cite{giannelli2023eg}. 
Finally, the pseudo-code of the reweighted least squares {adaptive spline} algorithm is reported in~\Cref{alg:rWLS-adaptive}.

\begin{algorithm2e}[!t]
	\SetAlgoLined
	\KwData{Point cloud $\{\bm x_i,\bm f_i\}_{i=1}^m$, the set of markers $K_I, K_{II}\subseteq\{1,\ldots,m\}$, the tolerances $\textup{tol}_I$, $\textup{tol}_{II}$ and $\epsilon > 0$, the penalization weight $\lambda\geq0$, a tensor product spline space $V = \spn\{\beta_1, \ldots, \beta_n\}$ and a maximum number hierarchical level $\mathrm{L}$;}
	\KwResult{$\bm{v} \in \mathcal{H}$ the reweighted adaptive least squares approximant.}
	Initialize the weights $\omega_i = 1$ and the point-wise errors $e_i = 1$ for each $i=1, \ldots, m$ and $\mathrm{loop}=0$\\
	\While{$\max_{i}{e_i} > \epsilon$ and $\mathrm{loop} < \mathrm{L}$}{%
		Solve the penalized weighted least squares problem \eqref{eq:wls-pen}
		\\
		compute the {pointwise} errors $e_i = \|\bm{v}(\bm x_i) - \bm f_i \|_2$;\\
		for each $i \in K_I$,
		if $e_i > \mathrm{tol}_I$, update $\omega_i = \omega_i \cdot \alpha(e_i)$,
		else $K_I = K_I \setminus \{i\}$;\\
		for each $i \in K_{II}$,
		if $e_i < \mathrm{tol}_{II}$, update $\omega_i = \omega_i \cdot \alpha(e_i)$,
		else $K_{II} = K_{II} \setminus \{i\}$;\\
		mark the domain elements where $e_i > \epsilon$;\\
		refine by dyadic split of the marked cells;\\
		update the data parameterization;\\
		set loop = loop$+$1\\
	}
	\KwRet{$\bm{v} \in \mathcal{H}$ which  solves the weighted adaptive spline least squares problem.}
	\caption{Reweighted adaptive least squares spline fitting.}
	\label{alg:rWLS-adaptive}
\end{algorithm2e}

\section{Numerical experiments}\label{sec:Numerical_Experiments}
In this section, we present a selection of numerical experiments to show the performance of the proposed fitting schemes. 
More precisely, in~\Cref{subsec:num_exp_rWLS-landmarks}, we show the effectiveness of the proposed reweighted least squares spline fitting method to recover the sharp features of the four different point sequences with type I markers, depicted in~\Cref{fig:landmarks}, and confront it with ordinary least squares spline fitting.
Moreover, in~\Cref{subsec:num_exp_rWLS-smoothing}, for the task of curve fitting, we compare the proposed method also with smoothing splines approximations.
In~\Cref{subsec:adaptive_rWLS},  we then extend the proposed method to adaptive spline spaces for surface reconstruction and suggest an automatic recognition of sharp features and type I markers.
Finally, in \Cref{subsec:adaptive_rWLS_nefertiti}, we show the performance of the adaptive algorithm when only type II markers are considered.

The univariate examples have been implemented in an 8-core laptop (Apple M2) with 8 GB RAM using \texttt{MATLAB R2023b}, specifically by employing the \texttt{Curve Fitting Toolbox}~\cite{curvefittingMatlab}. The adaptive surface approximations have been implemented within the open source  \texttt{C++} Geometry + Simulation (\texttt{G+Smo}) library~\cite{m2020}, by suitably extending the \texttt{gsHFitting} class.
Both libraries provide an efficient and robust implementation of least squares problems by exploiting suitable solvers for the linear systems at hand, which avoid the direct solution of a linear system of normal equations.

\subsection{Reweighted spline curve fitting with type I markers and comparison with ordinary least squares}
\label{subsec:num_exp_rWLS-landmarks}
To show the performance of~\Cref{alg:rWLS-landmarks}, we consider the points clouds and markers of type I depicted in~\Cref{fig:landmarks}.
The point clouds are similar to the ones presented in the following works: for ~\Cref{fig:landmarks}~(a)~\cite{ohrhallinger2021}, (b)~\cite{laube2020parametrization}, (c) and (d)~\cite{liu2015high}.
Our main goal is to improve the reconstruction while keeping the spline space intact.
Therefore,  each considered point sequence is parameterized accordingly to the \emph{uniform parameterization} and we compute the approximations of each dataset for the same spline space, i.e. with the same polynomial degree the same amount of interior nodes of the \emph{uniform} knot vector.	
We set the tolerance of the proposed method as $\textup{tol}_I=10^{-3}$ and compare the solution of the reweighted least squares problem (rWLS) with the solution of the ordinary least squares problem (LS) i.e., all the weights are equal to one.
The solutions are depicted in~\Cref{fig:solutions_landmarks}~(a),~(b),~(c), and (d), for the input data of ~\Cref{fig:landmarks}~(a),~(b),~(c),~and~(d) respectively. In particular, we can clearly see that the approximation using the reweighted least squares method (depicted in purple in \Cref{fig:solutions_landmarks}) shows a better approximation power with respect to the solution of the ordinary least squares problem (shown in red). 
In particular, one can see when zooming in that the sharp features (illustrated in black) are preserved during the curve fitting.
 Regarding the approximation error for each experiment, we report that, if considering \emph{all} the observation data, the rooted mean squares error (RMSE) and maximum error (MAX) are comparable with the ordinary least squares.
However, if we compare only how well we can preserve the sharp features, the rWLS method considerably outperforms the ordinary LS. More precisely, in~\Cref{table:MSE_section4}, we report the RMSE and MAX errors measured only for the set of type I marks, respectively for rWLS and LS.

\begin{figure}[!t]
	\captionsetup[subfigure]{labelformat=empty}
	\begin{center}
	\begin{subfigure}{0.45\textwidth}
		\centering
		\begin{tikzpicture}[zoomboxarray, black and white=cycle, zoomboxarray columns=2, zoomboxarray rows=1, zoomboxes below, connect zoomboxes]
			\node [image node]{\includegraphics[width=0.7\textwidth,height=0.15\textheight]{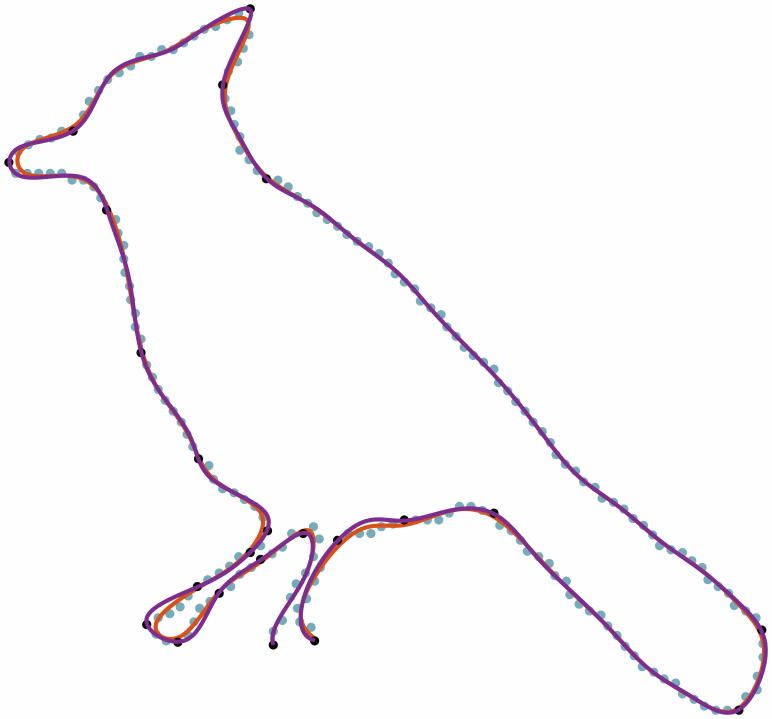} };
			\zoombox[magnification = 5,color code=black]{0.05,0.77}
			\zoombox[magnification = 5,color code=black]{0.30,0.95}
		\end{tikzpicture}
		\caption{(a) Degree 3 and $50$ interior knots}
	\end{subfigure}
	~
	\begin{subfigure}{0.45\textwidth}
		\centering
		\begin{tikzpicture}[zoomboxarray,black and white=cycle,  zoomboxarray columns=2, zoomboxarray rows=1, zoomboxes below,connect zoomboxes]
			\node [image node]{\includegraphics[width=0.7\textwidth, height= 0.15\textheight]{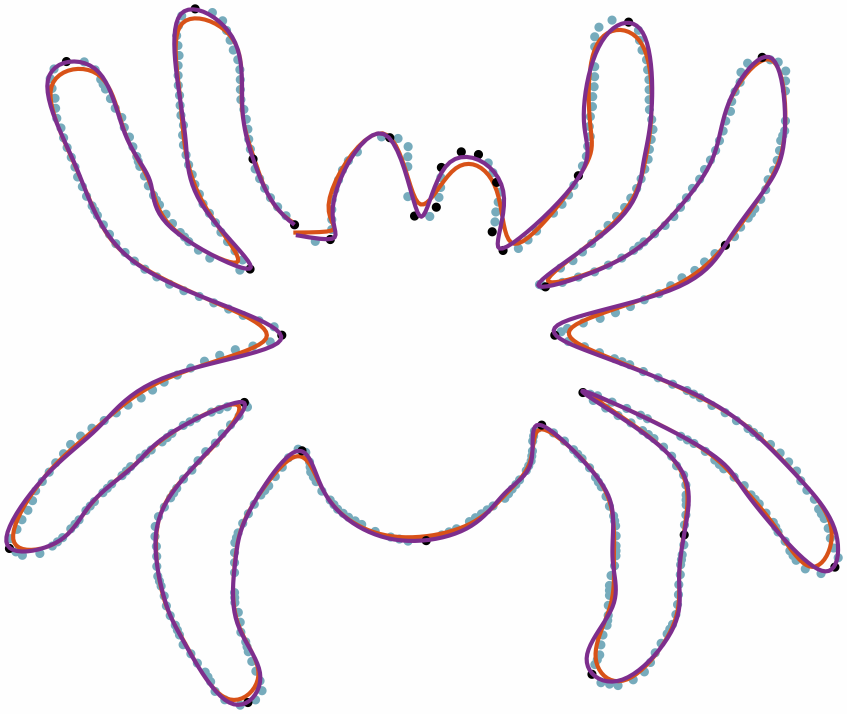} };
			\zoombox[magnification = 5,color code=black]{0.30,0.53}
			\zoombox[magnification = 5,color code=black]{0.55,0.72}
		\end{tikzpicture}
		\caption{(b) Degree 3 and $60$ interior knots}
	\end{subfigure}
	\\[10pt]
	\begin{subfigure}{0.45\textwidth}
		\centering
		\begin{tikzpicture}[zoomboxarray, black and white=cycle, zoomboxarray columns=1, zoomboxarray rows=1, zoomboxes below,connect zoomboxes]
			\node [image node]{\includegraphics[width=0.7\textwidth,height=0.15\textheight]{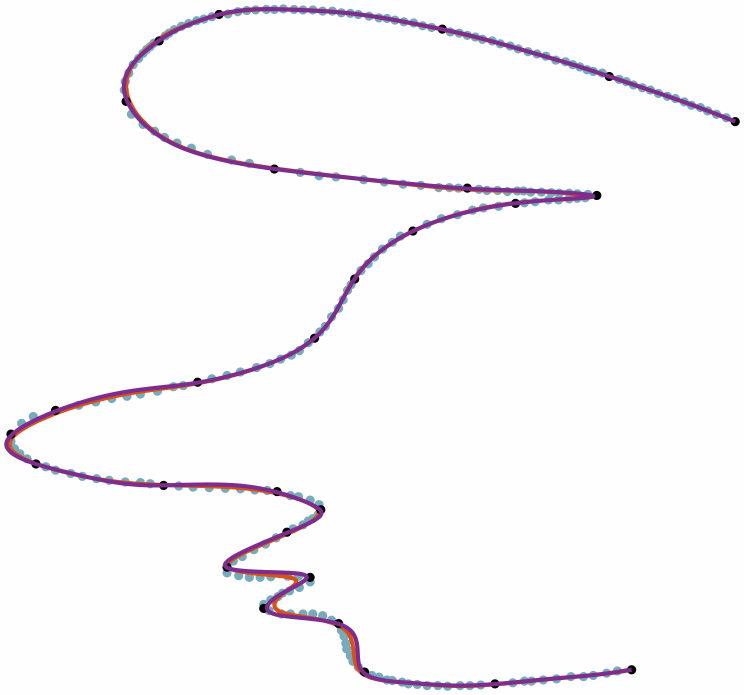} };
			\zoombox[magnification = 5,color code=black]{0.4,0.18}
		\end{tikzpicture}
		\caption{(c) Degree 3 and $40$ interior knots}
	\end{subfigure}
	~
	\begin{subfigure}{0.45\textwidth}
		\centering
		\begin{tikzpicture}[zoomboxarray, black and white=cycle, zoomboxarray columns=1, zoomboxarray rows=1, zoomboxes below,connect zoomboxes]
			\node [image node]{\includegraphics[width=0.7\textwidth,height=0.15\textheight]{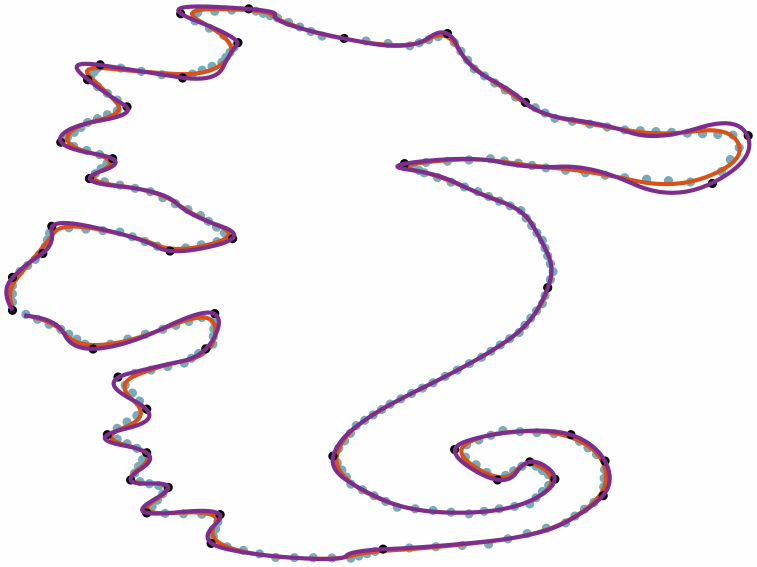} };
			\zoombox[magnification = 5,color code=black]{0.2,0.28}
		\end{tikzpicture}
		\caption{(d) Degree 2 and $84$ interior knots}
	\end{subfigure}
	\end{center}
	\caption{Curve fitting experiment described in~\Cref{subsec:num_exp_rWLS-landmarks} using ordinary LS (red) and rWLS (purple) for data with markers of type I.}
	\label{fig:solutions_landmarks}
\end{figure}

\begin{table}[!t]
	\sisetup{round-mode = places, round-precision = 3, scientific-notation = fixed, fixed-exponent = 0}
	\setlength\extrarowheight{2pt}
	\begin{center}
		\begin{tabular*}{\textwidth}{c@{\extracolsep{\fill}}ccccc}
			& Method & Iterations & Time & \multicolumn{1}{c}{$\textup{RMSE}$} & \multicolumn{1}{c}{$\textup{MAX}$} \\ \hline
			\multicolumn{1}{c|}{\multirow{2}{*}{(a)}}     & LS         & --   & -- &   $6.27 \times 10^{-3}$   &  $1.30\times 10^{-2}$    \\ 
			\multicolumn{1}{c|}{}                          & rWLS       & 1442   & \SI{1.347471}{\second} & $8.31 \times 10^{-4}$   &  $8.95\times 10^{-4}$   \\ \hline
			\multicolumn{1}{c|}{\multirow{2}{*}{(b)}}   & LS         & --    & -- &  $2.25\times 10^{-3}$   & $3.94\times 10^{-3}$    \\ 
			\multicolumn{1}{l|}{}                          & rWLS       & 1250   & \SI{1.786006}{\second} & $6.63\times 10^{-4}$   & $1.78\times 10^{-3}$     \\ \hline
			\multicolumn{1}{c|}{\multirow{2}{*}{(c)}}     & LS        & --    & -- &  $4.63\times 10^{-3}$   &  $1.22\times 10^{-2}$    \\ 
			\multicolumn{1}{l|}{}                          & rWLS         & 30   & \SI{0.033928}{\second}& $1.47\times 10^{-3}$  &  $3.54\times 10^{-3}$     \\ \hline
			\multicolumn{1}{c|}{\multirow{2}{*}{(d)}} & LS         & --  & --  &  $5.94\times 10^{-3}$   &  $1.25\times 10^{-2}$    \\ 
			\multicolumn{1}{l|}{}                          & rWLS         & 81 & \SI{0.137070}{\second}  &  $3.77\times 10^{-4}$  &  $8.81\times 10^{-4}$   
		\end{tabular*}
		\caption{RMSE and MAX errors for type I markers of LS and rWLS solutions of the experiment described in~\Cref{subsec:num_exp_rWLS-landmarks}. For rWLS also the number of iterations (Iterations) is reported.}
		\label{table:MSE_section4}
	\end{center}
\end{table}

\subsection{Reweighted spline curve fitting with type I markers and comparison with smoothing splines}
\label{subsec:num_exp_rWLS-smoothing}

In addition to ordinary least squares, another widely used method for curve fitting is smoothing splines \cite{gu1992, gu2013}. This method is derived as the solution to the penalized least squares problem \eqref{eq:wls-pen}, in terms of cubic splines with knots corresponding to the $x$-coordinates of the observations. Therefore, we considered three additional point clouds sampled from the following functions:
\begin{equation*}
\small
f_1(x)  = \left|\frac{9\sin(3\pi x)}{\tanh(-1.5x+1)+1}\right|, \;
f_2(x)   = \frac{1}{0.02\sqrt{\pi}} \exp\left\{-\left(\frac{x-0.5}{0.02}\right)^2 \right\}, \;
f_3(x) = \tanh\left(\frac{\cos(2\pi x)}{0.05}\right)
\end{equation*}
and create $3$ point clouds with $62,88,71$ points, respectively.
Specifically, we consider an irregular distribution of the abscissas to obtain more observations around the functions' sharp features.

We compute the smoothing spline fitting using the \texttt{fit} function in \texttt{MATLAB} and use the default ``interesting range'' value of the smoothing weight $\lambda$, which depends on the abscissa distribution of each experiment. Moreover, we set the weights of \eqref{eq:wls-pen} to $1$.
To ensure a fair comparison, for the rWLS algorithm, we fixed the spline degree to $3$ and selected the knot vectors following the ``averaging'' technique described in~\cite[Eq.~(9.8)]{Piegl1997}, known in the literature also as NKTP, with $37$, $47$, and $47$ interior knots, respectively. This technique has been chosen since it also considers the values of the abscissas as in the case of smoothing splines. 

The solutions obtained using rWLS from \Cref{alg:rWLS-landmarks} are depicted in \Cref{fig:smoothing_example}. Additionally, we report the rooted mean squared error (RMSE) and maximum error (MAX) measured only for the set of type I marks for both rWLS and smoothing splines in \Cref{table:MSE_smoothing}. To summarise, in the presence of sharp features, better results can be obtained using the proposed rWLS from \Cref{alg:rWLS-landmarks} with fewer degrees of freedom compared to the smoothing spline technique. 
}

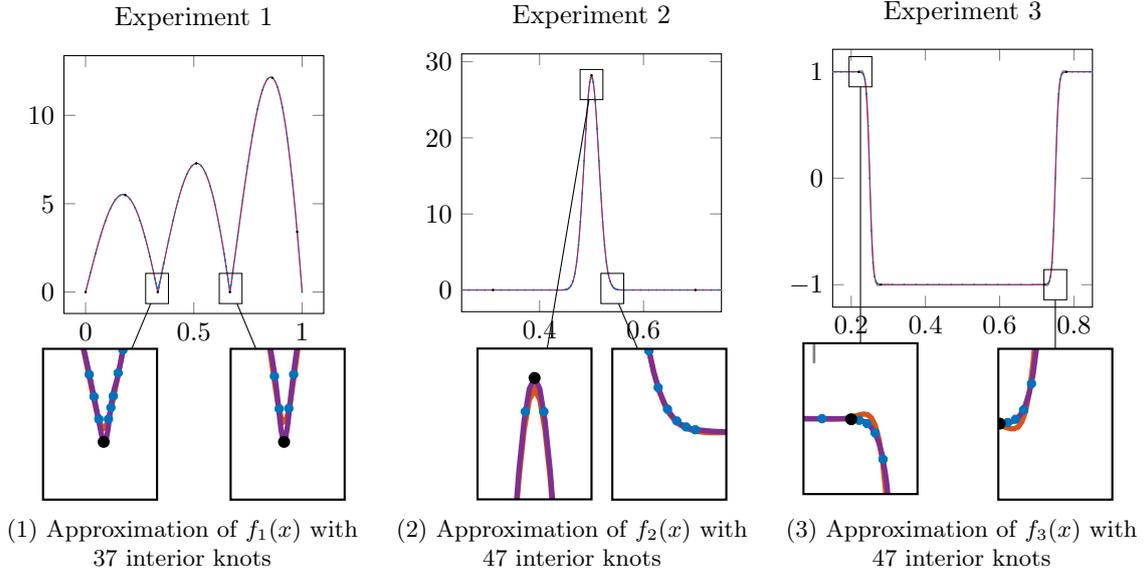
\begin{figure}[!t]
\renewcommand\thesubfigure{\arabic{subfigure}}
\captionsetup[subfigure]{justification=centering}
\begin{center}
\begin{subfigure}{0.3\textwidth}
\begin{center}
\begin{tikzpicture}
\begin{scope}[spy using outlines={magnification=5, width=1.5cm,height=2cm,connect spies}]
\begin{axis}[title= Experiment 1, legend style={draw = none}, legend pos =  north west, width =5cm, height = 5cm ]
\addplot[only marks, color = pointCloud, mark = *, mark size = 0.15pt] table [x=x,y=y, col sep=comma]{pointCloud_example_4.txt};
\addplot[color = theirSolution] table [x=x,y=f_fit, col sep=comma]{solution_example_4.txt};
\addplot[color = ourSolution] table [x=x,y=f_rwls, col sep=comma]{solution_example_4.txt};
\addplot[only marks, color = black, mark = *, mark size = 0.25pt] table [x=x,y=y, col sep=comma]{landmarks_example_4.txt};
\path (0.33, 0.2) coordinate (X);
\path (0.67, 0.2) coordinate (Y);
\end{axis}
 \spy [black and white=cycle] on (X) in node (zoom) [left] at ([xshift=0cm,yshift=-1.8cm]X);
 \spy [black and white=cycle] on (Y) in node (zoom) [right] at ([xshift=0cm,yshift=-1.8cm]Y);
 \end{scope}
\end{tikzpicture}
\caption{Approximation of $f_1(x)$ with $37$ interior knots}
\end{center}
\end{subfigure}
~
\begin{subfigure}{0.3\textwidth}
\begin{center}
\begin{tikzpicture}
\begin{scope}[spy using outlines={magnification=5, width=1.5cm,height=2cm,connect spies}]
\begin{axis}[title= Experiment 2, legend style={draw = none}, legend pos =  north west, width =5cm , height= 5cm, xmin = 0.25, xmax = 0.75]
\addplot[only marks, color = pointCloud, mark = *, mark size = 0.15pt] table [x=x,y=y, col sep=comma]{pointCloud_example_10.txt};
\addplot[color = theirSolution] table [x=x,y=f_fit, col sep=comma]{solution_example_10.txt};
\addplot[color = ourSolution] table [x=x,y=f_rwls, col sep=comma]{solution_example_10.txt};
\addplot[only marks, color = black, mark = *, mark size = 0.25pt] table [x=x,y=y, col sep=comma]{landmarks_example_10.txt};
\path (0.5, 27) coordinate (X);
\path (0.54, 0.2) coordinate (Y);
\end{axis}
 \spy [black and white=cycle] on (X) in node (zoom) [left] at ([xshift=0cm,yshift=-4.5cm]X);
 \spy [black and white=cycle] on (Y) in node (zoom) [right] at ([xshift=0cm,yshift=-1.8cm]Y);
 \end{scope}
\end{tikzpicture}
\caption{Approximation of $f_2(x)$  with $47$ interior knots}
\end{center}
\end{subfigure}
~
\begin{subfigure}{0.3\textwidth}
\begin{tikzpicture}
\begin{scope}[spy using outlines={magnification=5, width=1.5cm,height=2cm,connect spies}]
\begin{axis}[title= Experiment 3, legend style={draw = none}, legend pos =  north west, width =5cm , height= 5cm, xmin = 0.15, xmax = 0.85]
\addplot[only marks, color = pointCloud, mark = *, mark size = 0.15pt] table [x=x,y=y, col sep=comma]{pointCloud_example_12.txt};
\addplot[color = theirSolution] table [x=x,y=f_fit, col sep=comma]{solution_example_12.txt};
\addplot[color = ourSolution] table [x=x,y=f_rwls, col sep=comma]{solution_example_12.txt};
\addplot[only marks, color = black, mark = *, mark size = 0.25pt] table [x=x,y=y, col sep=comma]{landmarks_example_12.txt};
\path (0.225, 1) coordinate (X);
\path (0.75, -1) coordinate (Y);
\end{axis}
 \spy [black and white=cycle] on (X) in node (zoom) [below] at ([xshift=0cm,yshift=-3.6cm]X);
 \spy [black and white=cycle] on (Y) in node (zoom) [below] at ([xshift=0cm,yshift=-0.85cm]Y);
\end{scope}
\end{tikzpicture}
\caption{Approximation of $f_3(x)$  with $47$ interior knots}
\end{subfigure}
\caption{Curve fitting experiment described in~\Cref{subsec:num_exp_rWLS-smoothing} using smoothing spline (red) and rWLS (purple) for data with markers of type I (black dots).}
\label{fig:smoothing_example}
\end{center}
\end{figure}

\begin{table}[!t]
	\sisetup{round-mode = places, round-precision = 3, scientific-notation = fixed, fixed-exponent = 0}
	\setlength\extrarowheight{2pt}
	\begin{center}
		\begin{tabular*}{\textwidth}{c@{\extracolsep{\fill}}ccccc}
			& Method & Iterations & Time &  \multicolumn{1}{c}{$\textup{RMSE}$} & \multicolumn{1}{c}{$\textup{MAX}$} \\ \hline
			\multicolumn{1}{c|}{\multirow{2}{*}{(1)}}     & smoothing spline         & --    & \SI{0.20188}{\second}& $1.02 \times 10^{-1}$   &  $2.21\times 10^{-1}$    \\ 
			\multicolumn{1}{c|}{}                          & rWLS       & 10 & \SI{0.086209}{\second}  &  $9.13 \times 10^{-6}$   &  $1.86\times 10^{-5}$   \\ \hline
			\multicolumn{1}{c|}{\multirow{2}{*}{(2)}}   & smoothing spline         & --  & \SI{0.18782}{\second}    &   $1.36\times 10^{-1}$   & $3.04\times 10^{-1}$    \\ 
			\multicolumn{1}{l|}{}                          & rWLS       & 7  & \SI{0.027153}{\second}   &  $1.43\times 10^{-5}$   & $3.18\times 10^{-5}$     \\ \hline
			\multicolumn{1}{c|}{\multirow{2}{*}{(3)}}     & smoothing spline        & --  & \SI{0.14277}{\second}  &  $2.96\times 10^{-3}$   &  $3.86\times 10^{-3}$    \\ 
			\multicolumn{1}{l|}{}                          & rWLS         & 4  & \SI{0.010721}{\second}   &  $1.30\times 10^{-5}$  &  $2.58\times 10^{-5}$     \\ \hline
		\end{tabular*}
		\caption{RMSE and MAX errors for type I markers of smoothing spline and rWLS solutions of the experiment described in~\Cref{subsec:num_exp_rWLS-smoothing}. For rWLS also the number of iterations (Iterations) is reported.}
		\label{table:MSE_smoothing}
	\end{center}
\end{table}

\subsection{Reweighted adaptive spline fitting for type I markers}
\label{subsec:adaptive_rWLS}
\begin{figure}[!t]
	\centering
	\resizebox{\linewidth}{!}
	{
		\begin{tikzpicture}
			\begin{axis}[width=0.45\linewidth, name=plot1, xlabel={DOFs}, ylabel = {MAX}, ymode=log, legend style={draw=none}]
				\addplot+[color=red, thick, mark=square, mark size= 1pt, mark options={solid}] table [x = dofs, y=max, col sep=comma] {LS_results.csv};
				\addlegendentry{LS}
				\addplot+[color=purple, thick, mark=triangle, mark size= 1pt, mark options={solid}] table [x = dofs, y=max, col sep=comma] {IR_results.csv};
				\addlegendentry{rWLS}
			\end{axis}
			\node (wls-mesh) [right = 2pt of plot1] {
					\includegraphics[width=0.28\linewidth, trim = 200 150 200 150, clip]{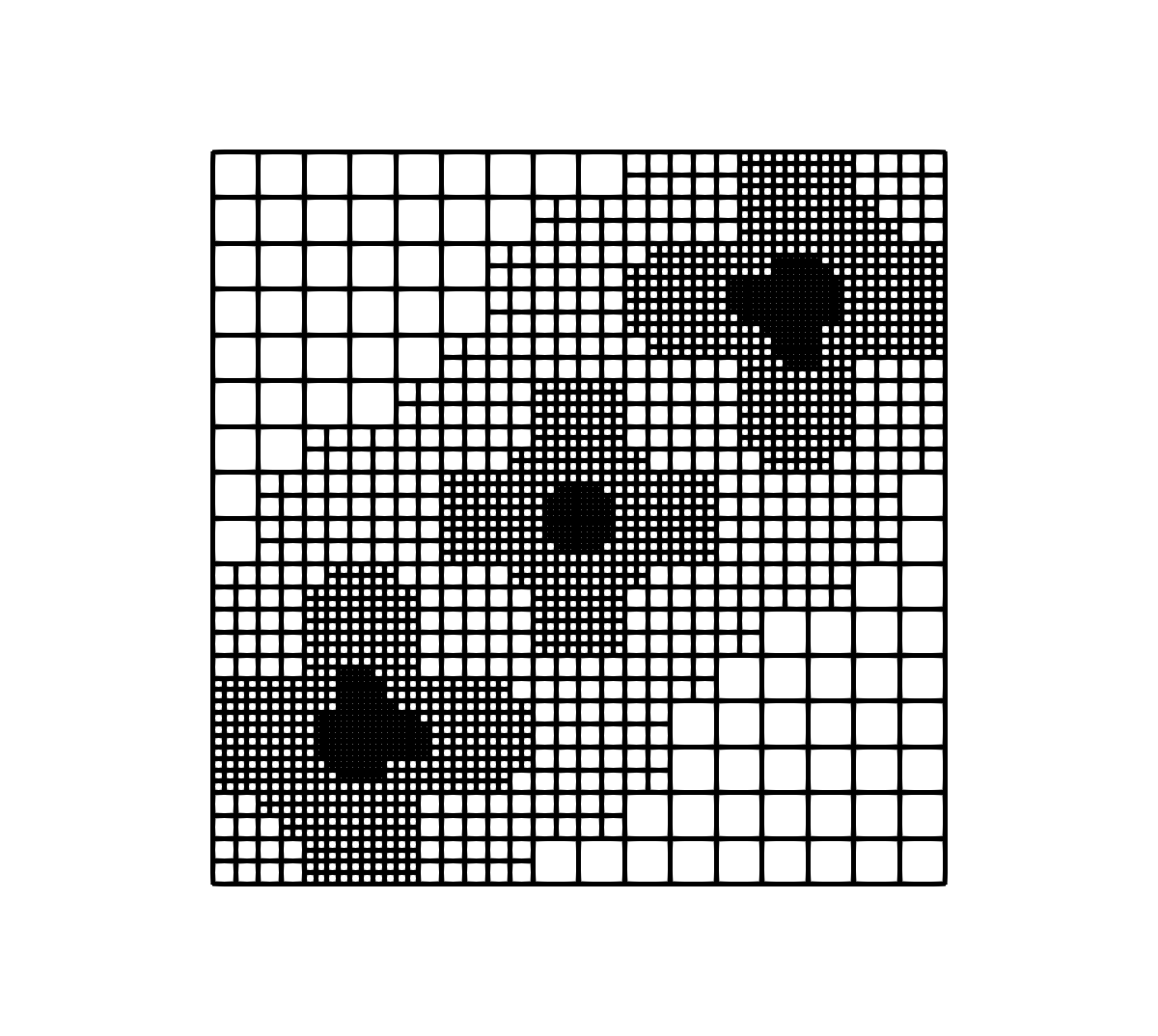}
				};
			\node (wls-geo) [right = 0pt of wls-mesh]{
				\includegraphics[width=0.3\linewidth, trim = 250 100 200 200, clip]{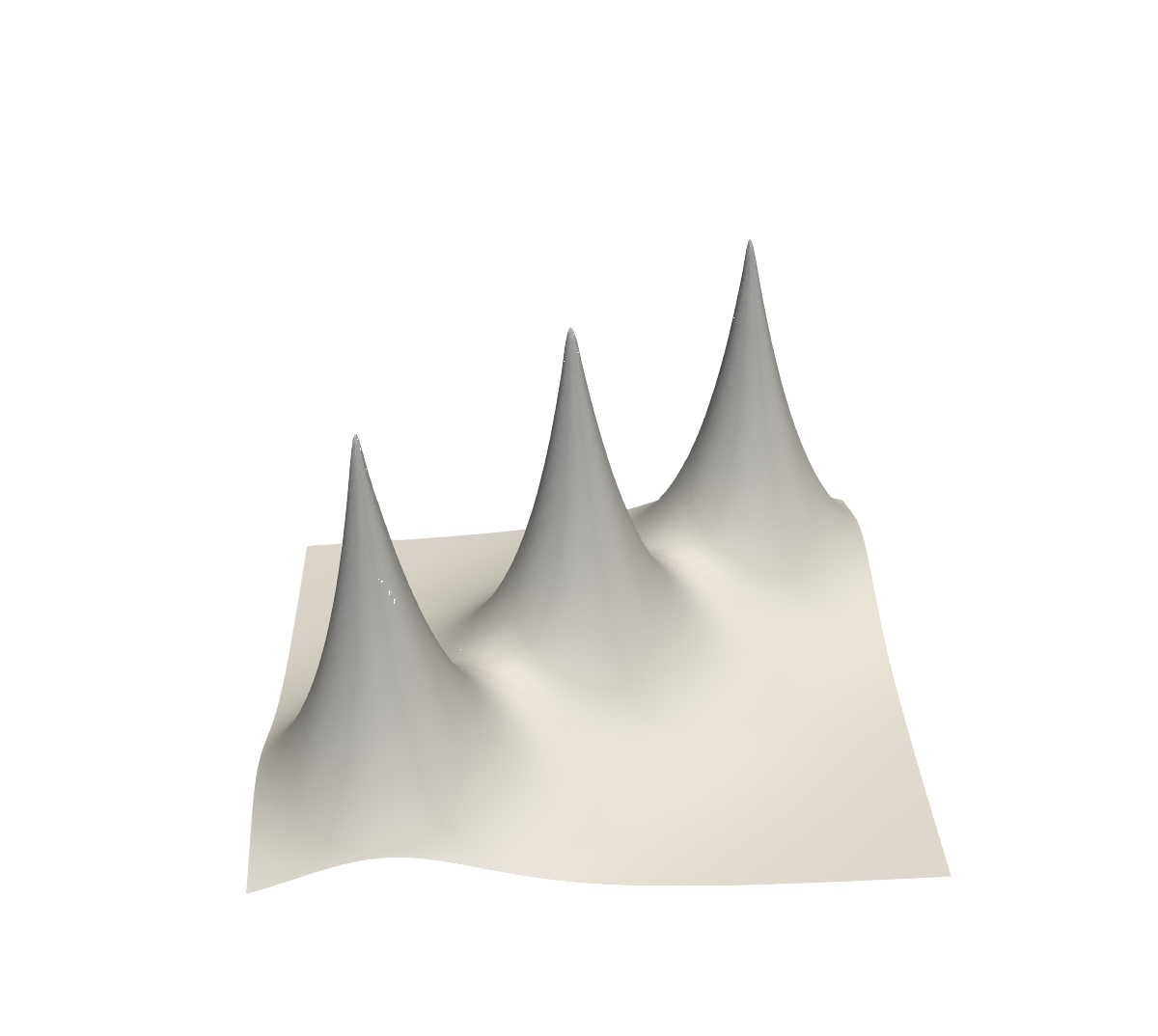}
				};
		\end{tikzpicture}
	}
	\caption{MAX error w.r.t. DOFs for the ordinary LS and rWLS fitting scheme THB-splines (left). Hierarchical mesh (center) and geometry (right) are obtained in output from the rWLS method.}
	\label{fig:exm-adaptive-3peaks}
\end{figure}

In this numerical experiment we apply \Cref{alg:rWLS-adaptive}, combined with truncated hierarchical B-splines (THB-splines) as in \cite{giannelli2023eg}, to a point cloud obtained by sampling $100 \times 100$ gridded data from function \eqref{eq:3peaks} used in \Cref{exm:box-splines}. 
This function is characterized by sharp features which we want the final model to capture, i.e. we would deal with type I markers only, hence $K_{II} = \emptyset$. 
In particular to initialize the set $K_I$ we perform an ordinary least squares fitting with tensor-product B-splines, of bi-degree 3 and a $15\times 15$ mesh, and individuate the data sites whose approximation error is above a certain threshold $\epsilon$, namely $K_I \coloneqq \{ i \in \{1, \ldots, m\}: \left\| v(\bm x_i) - f(\bm x_i)\right\|_2 > \epsilon\}$. 
Given the data points, the markers $K_I$, the tolerances and the initial tensor-product B-spline space, we perform \Cref{alg:rWLS-adaptive} with $\mathrm{tol}_I = 10\epsilon$ augmenting the weights value by 25\% at each iteration of the adaptive loop. 
We then compare the output of the reweighted adaptive least squares fitting method (rWLS) with the result of the ordinary adaptive least squares fitting (LS), namely by keeping all the weights equal to one all along the iterative procedure.
More precisely, we compare the maximum error (MAX) of the final approximations with respect to the number of degrees of freedom (DOFs), choosing $\epsilon$ to achieve a comparable number of DOFs. 
The comparison is reported in \Cref{fig:exm-adaptive-3peaks} (left), together with the final hierarchical mesh (center) and THB-spline approximation (right) obtained with the rWLS method. 
More precisely, LS and rWLS methods have the same starting point and as long as the adaptive refinement proceeds, we can notice that rWLS can register a smaller MAX error with fewer DOFs in comparison to LS.

\subsection{Reweighted adaptive spline fitting for type II markers}
\label{subsec:adaptive_rWLS_nefertiti}
\begin{figure}[!t]
	\captionsetup[subfigure]{labelformat=empty}
	\centering
	\hspace*{-9cm}\resizebox{0.35\linewidth}{!}{
		\begin{subfigure}{0.5\linewidth}
			\centering
			\begin{tikzpicture}[zoomboxarray,black and white=cycle,  zoomboxarray columns=1.7, zoomboxarray rows=1.7, zoomboxes right, connect zoomboxes]
				\node (pc) {\includegraphics[width=\linewidth,trim= 600 400 600 700,clip]{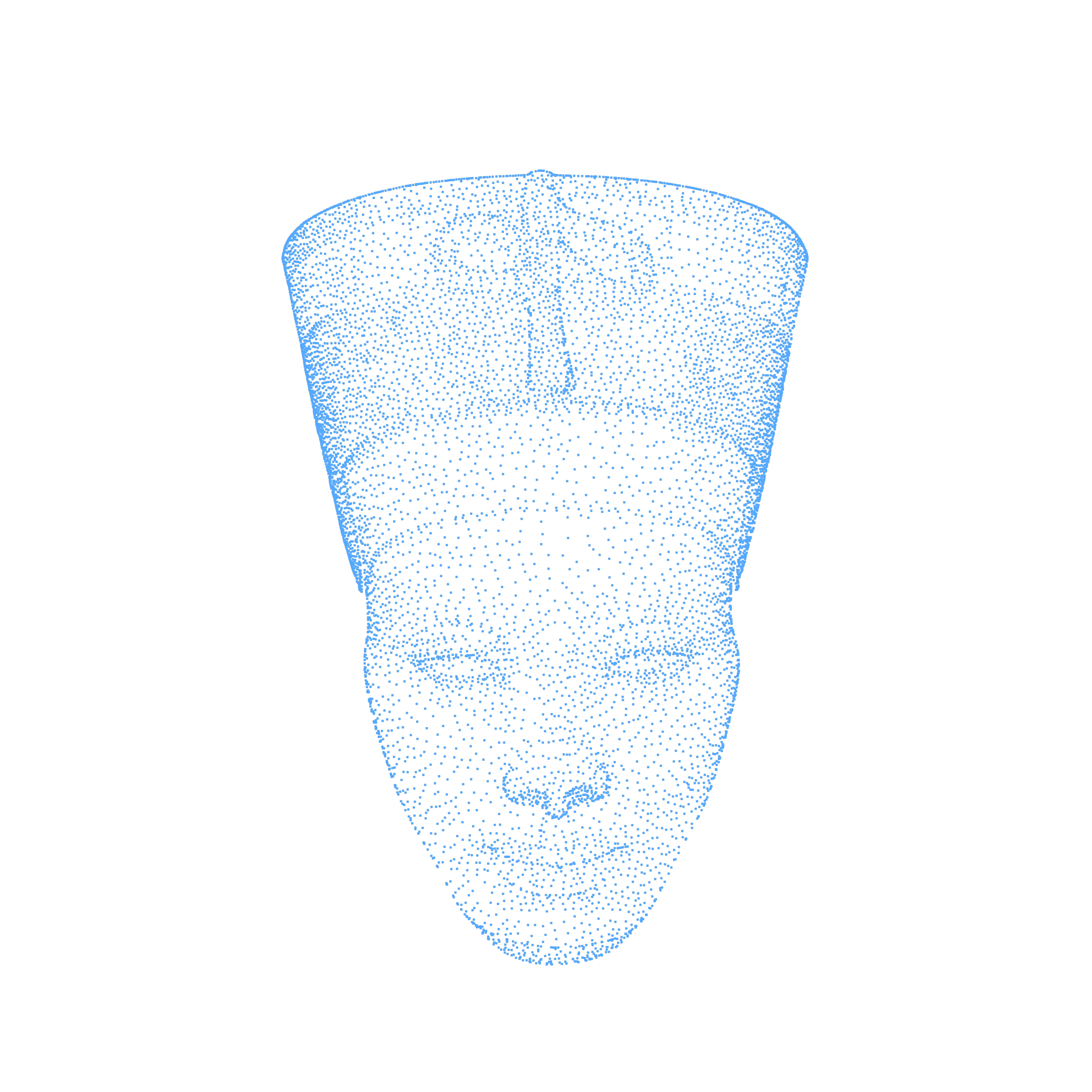}};
				\node [image node, right = 0pt of pc] {\includegraphics[width=\linewidth,trim= 600 400 600 700,clip]{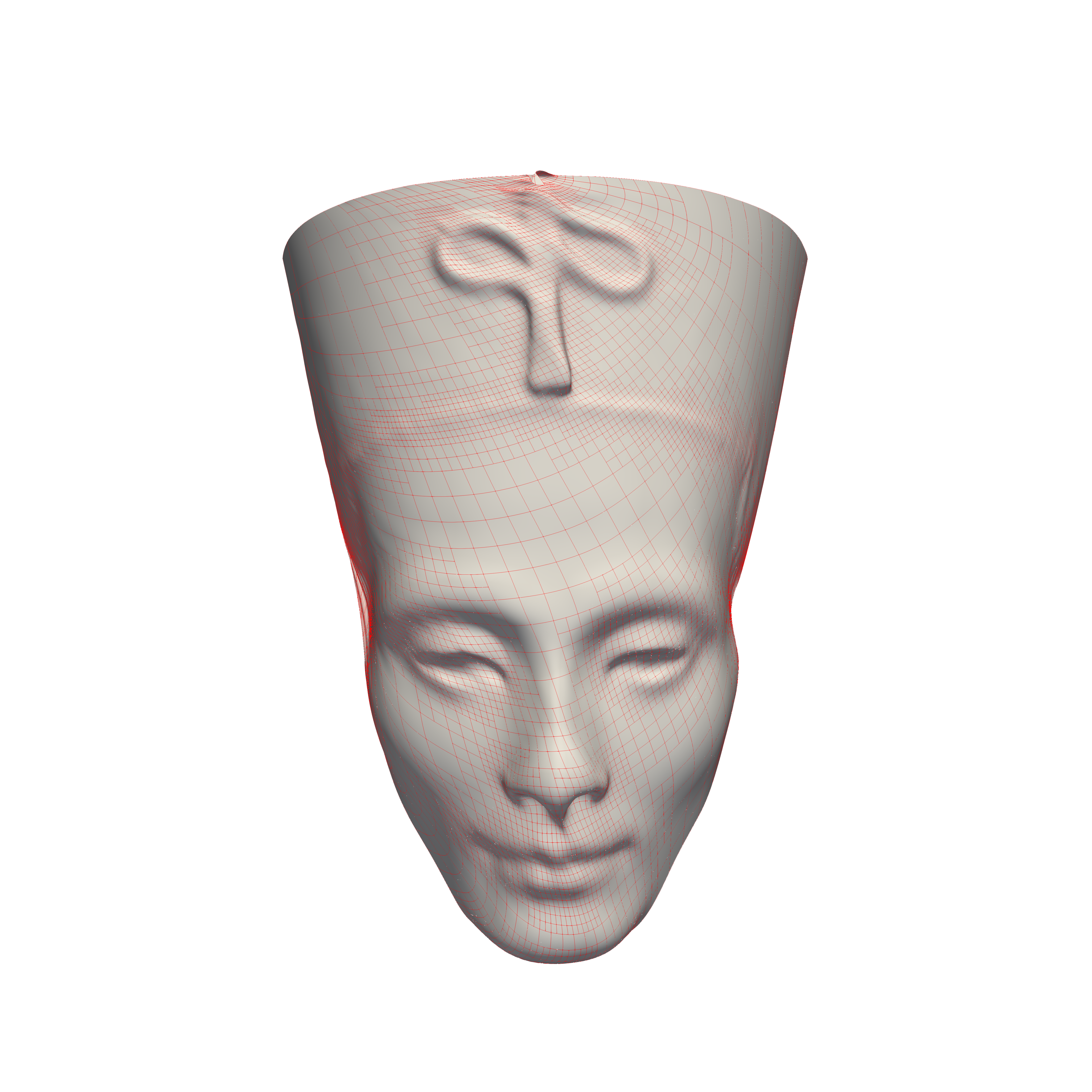}};
				\zoombox[magnification = 2,color code=black]{1.0125,0.345}
			\end{tikzpicture}
		\end{subfigure}
	}

		\hspace*{-9cm}\resizebox{0.35\linewidth}{!}{
		\begin{subfigure}{0.5\linewidth}
			\centering
		\begin{tikzpicture}[zoomboxarray,black and white=cycle,  zoomboxarray columns=1.7, zoomboxarray rows=1.7, zoomboxes right, connect zoomboxes]
			\node (pcw) {\includegraphics[width=\linewidth,trim= 600 400 600 700,clip]{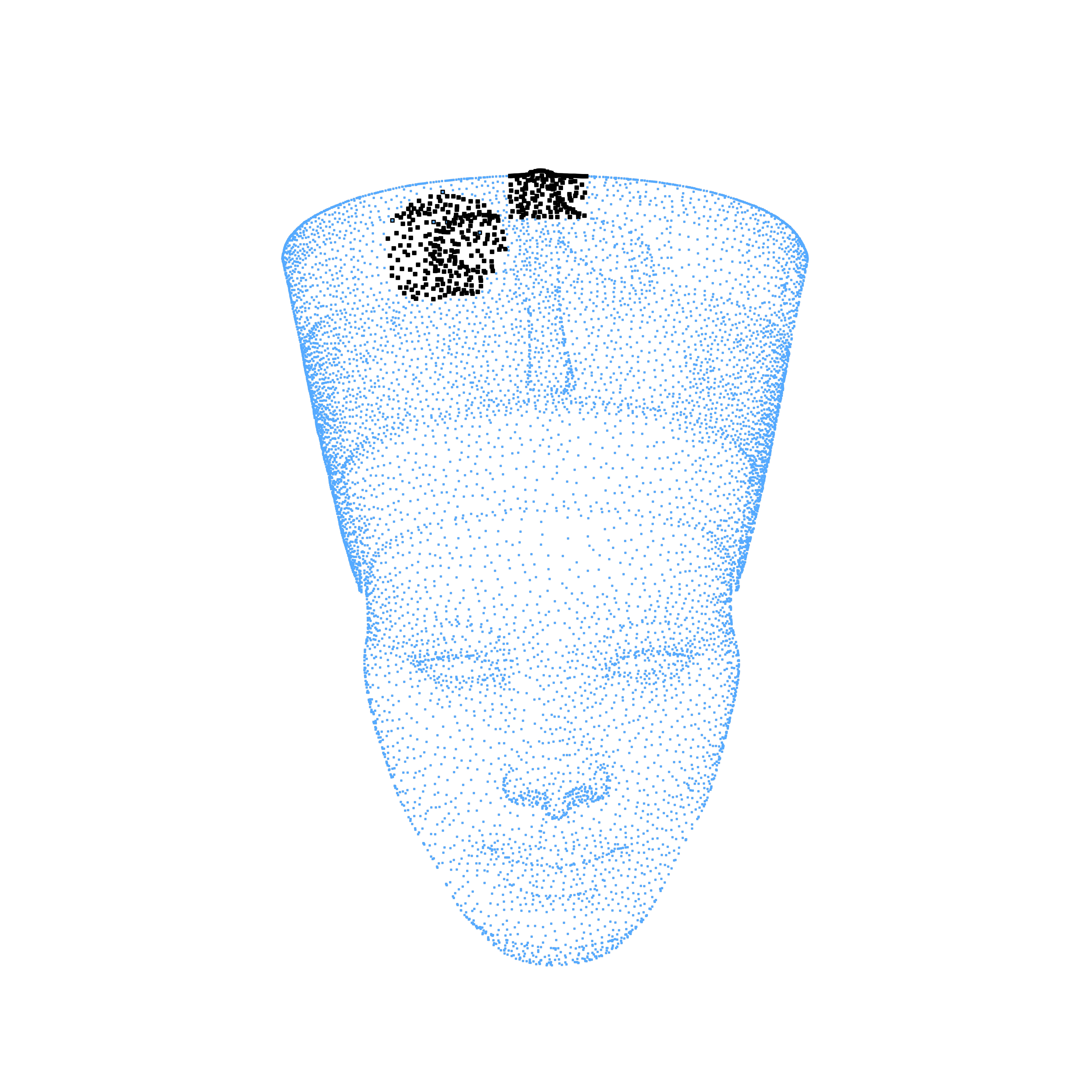}};
			\node [image node, right = 0pt of pcw]{\includegraphics[width=\linewidth,trim= 600 400 600 700,clip]{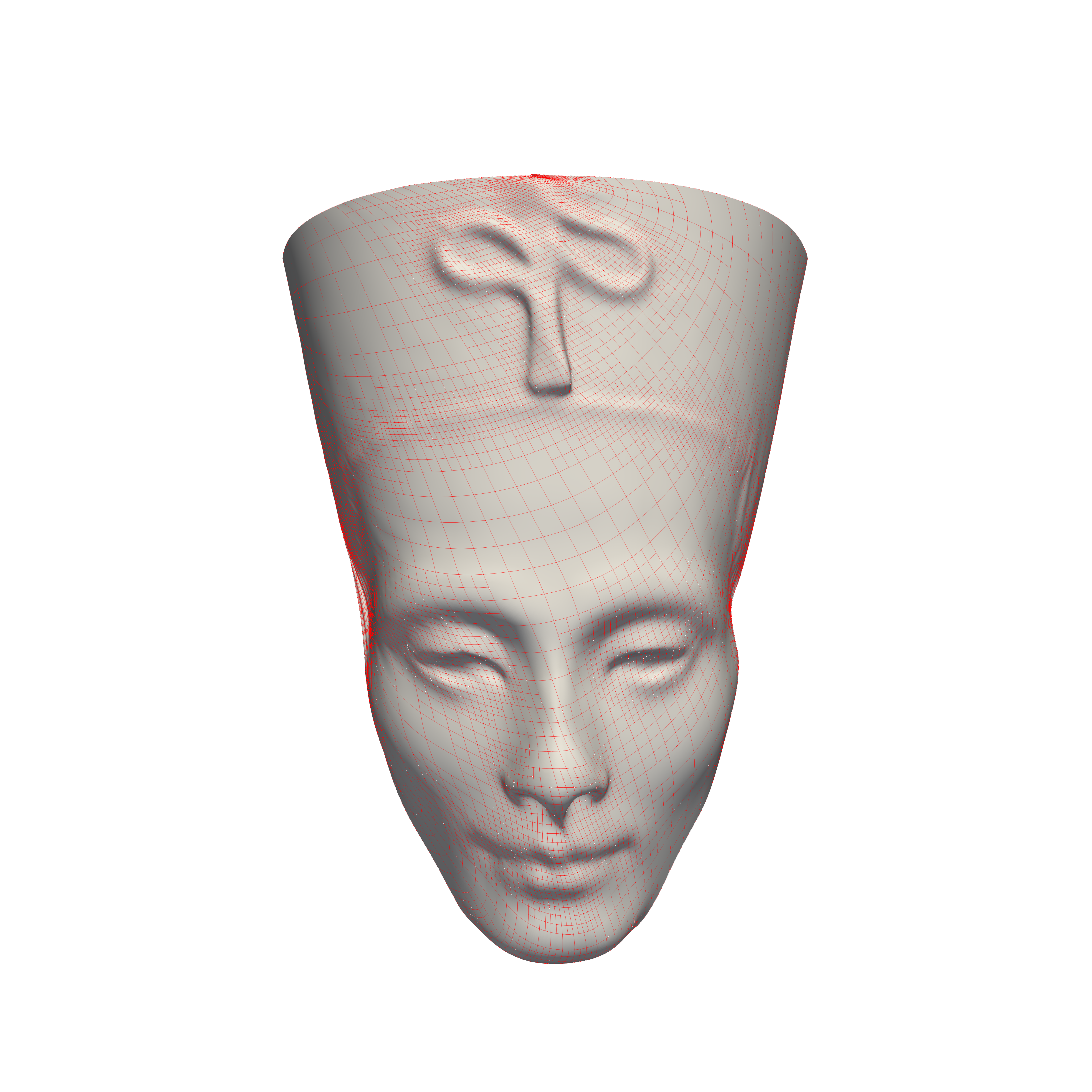}};
			\zoombox[magnification = 2,color code=black]{1.0125,0.345}
		\end{tikzpicture}
		\end{subfigure}
	}
		\caption{Adaptive surface fitting experiment described in~\Cref{subsec:adaptive_rWLS_nefertiti} using adaptive LS (top) and rWLS (bottom) for data with markers of type II. The data point clouds and the marked points are also shown.}\label{fig:weights-ii}
		\end{figure}

In this final numerical experiment, we perform the reweighted adaptive spline least squares fitting scheme to approximate a \emph{scattered} point cloud of 9000 points representing a Nefertiti bust, illustrated in \Cref{fig:weights-ii} (top-left). 
In particular, these data require a parameterization and some of them result in being corrupted. 
To address the first problem, we employ the data-driven parameterization method developed in \cite{giannelli2023nc}. 
To tackle the presence of corrupted data, we exploit the use of the weights associated with the data. 
More precisely, together with the points, we are given a set of markers of the second type $K_{II}$, highlighted in black in \Cref{fig:weights-ii} (bottom-left) which we use to perform the reweighted adaptive scheme. 
Moreover, due to the complexity of the problem, we decided not to update the markers $K_{II}$ but to keep them fixed along the iterative scheme, by setting $\mathrm{tol}_{II}$ sufficiently small. 
We then compare the results obtained by the proposed reweighted adaptive least squares scheme (rWLS) with the ordinary adaptive least squares scheme (LS), namely keeping all the weights set to one. 
Such a comparison is illustrated in \Cref{fig:weights-ii}, where we show on the top the hierarchical approximation obtained by LS and on the bottom the one obtained by rWLS. 
In particular, the region corresponding to the corrupted data is better approximated by the rWLS method, which results free of self-intersections and it is {smoother} than the one obtained with LS. 
Finally, the MAX error registered by LS is  3.28e$-$2, whereas the MAX of rWLS is 1.57e$-$2.

\section{Conclusions}
We provided new insights into the reweighted least squares method for the approximation of a given set of observations. In particular, our theoretical study revealed that any reweighted least squares model can be interpreted as a convex combination of suitable interpolation models.
Furthermore, we exploited the weights for spline fitting problems, such as adaptive spline constructions. 
Thereby, we provided a strategy to automatically update the weights within the fitting scheme, either to emphasize data marked as sharp features or to smoothen data marked as corrupted. 
In the present work,
we assumed that the given data have been previously processed, hence their features or the presence of noise have already been identified. This information was exploited in the fitting schemes via an iterative update of the weights.
In future research, it is of interest to automatically classify the data marker used in the spline approximation scheme
in order to drive the model design process.
Another interesting future research direction is the application of the proposed method to a non-geometric modelling framework, e.g., the approximation of time-series data. In addition, in this context, it would also be interesting to investigate how the choice of the smoothing parameter affects the goodness of the fit from a statistical point of view.

\paragraph{\textbf{CRediT authorship contribution statement}}
\textbf{Carlotta Giannelli:} conceptualization, methodology, investigation, formal analysis, writing -- original draft, review and editing, resources.
\textbf{Sofia Imperatore:} conceptualization, methodology, software, validation, investigation, writing -- original draft, review and editing, visualization.
\textbf{Lisa Maria Kreusser:} conceptualization, methodology, writing -- original draft, review and editing, formal analysis, resources.
\textbf{Estefan{\'i}a Loayza-Romero:} conceptualization, methodology, software, writing -- original draft, review and editing, visualization.
\textbf{Fatemeh Mohammadi:} conceptualization, resources.
\textbf{Nelly Villamizar:} writing -- original draft, reviewing and editing.

\paragraph{\textbf{Acknowledgments}}
The authors would like to acknowledge the support provided by the 4th WiSh: Women in Shape Analysis Research Workshop. 
This collaboration began during the workshop, and we are deeply grateful for the opportunity to work with fellow researchers in the field.
CG and SI are members of the INdAM group GNCS, whose support is gratefully acknowledged. LMK acknowledges support from Magdalene College, Cambridge (Nevile Research Fellowship). ELR work was partially funded by the Deutsche Forschungsgemeinschaft (DFG, German Research Foundation) under Germany's Excellence Strategy EXC 2044 –390685587, Mathematics Münster: Dynamics–Geometry–Structure. FM was partially supported by the FWO grants (G0F5921N, G023721N), the KU Leuven iBOF/23/064 grant, and the UiT Aurora MASCOT project. NV was supported by the UK Engineering and Physical Sciences Research Council (EPSRC) New Investigator Award EP/V012835/1.

\bibliographystyle{abbrv}
\bibliography{sample}

\end{document}